%% file: MarchildonZingg_CstrBayeOptz.tex
\documentclass[AMA,Times1COL]{WileyNJDv5} %STIX1COL,STIX2COL,STIXSMALL

\articletype{Research article}%

\received{Date Month Year}
\revised{Date Month Year}
\accepted{Date Month Year}
\journal{Journal}
\volume{00}
\copyyear{2025}
\startpage{1}

\input{NewCommands}

\begin{document}

\title{A Framework for Nonlinearly-Constrained Gradient-Enhanced Local Bayesian Optimization with Comparisons to Quasi-Newton Optimizers}

\author[1,2]{Andr\'e L.\ Marchildon}

\author[1]{David W.\ Zingg}

\authormark{Marchildon and Zingg}
\titlemark{Nonlinearly-Constrained Gradient-Enhanced Local Bayesian Optimization}

\address[1]{\orgdiv{Institute for Aerospace Studies}, \orgname{University of Toronto}, \orgaddress{\state{Ontario}, \country{Canada}}}

\address[2]{This work was completed prior to the author beginning his employment at Amazon}

\corres{Corresponding author Andr\'e L.\ Marchildon: \email{andre.marchildon@alumni.utoronto.ca}}

\fundingInfo{Natural Sciences and Engineering Research Council of Canada and the Ontario Graduate Scholarship Program.}

\abstract[Abstract]{Bayesian optimization is a popular and versatile approach that is well suited to solve challenging optimization problems. Their popularity comes from their effective minimization of expensive function evaluations, their capability to leverage gradients, and their efficient use of noisy data. Bayesian optimizers have commonly been applied to global unconstrained problems, with limited development for many other classes of problems. In this paper, two alternative methods are developed that enable rapid and deep convergence of nonlinearly-constrained local optimization problems using a Bayesian optimizer. The first method uses an exact augmented Lagrangian and the second augments the minimization of the acquisition function to contain additional constraints. Both of these methods can be applied to nonlinear equality constraints, unlike most previous methods developed for constrained Bayesian optimizers. The new methods are applied with a gradient-enhanced Bayesian optimizer and enable deeper convergence for three  nonlinearly-constrained unimodal optimization problems than previously developed methods for constrained Bayesian optimization. In addition, both new methods enable the Bayesian optimizer to reach a desired tolerance with fewer function evaluations than popular quasi-Newton optimizers from SciPy and MATLAB for unimodal problems with 2 to 30 variables. The Bayesian optimizer had similar results using both methods. It is recommended that users first try using the second method, which adds constraints to the acquisition function minimization, since its parameters are more intuitive to tune for new problems
\footnote{Content from this paper was previously published in this thesis: Marchildon, A. L. Aug. 2024. “The Development of a Versatile and Efficient Gradient-Enhanced Bayesian Optimizer for Nonlinearly Constrained Optimization with Application to Aerodynamic Shape Optimization”. PhD Thesis. Toronto, Canada: University of Toronto. URL: https://utoronto.scholaris.ca/items/2037ae09-d3f7-4781-a8b6-71c2743d0ef8.}.}


%\jnlcitation{\cname{%
%\author{Taylor M.},
%\author{Lauritzen P},
%\author{Erath C}, and
%\author{Mittal R}}.
%\ctitle{On simplifying ‘incremental remap’-based transport schemes.} \cjournal{\it J Comput Phys.} \cvol{2021;00(00):1--18}.}
 
\maketitle

\renewcommand\thefootnote{}
%\footnotetext{\textbf{Abbreviations:} ANA, anti-nuclear antibodies; APC, antigen-presenting cells; IRF, interferon regulatory factor.}

\renewcommand\thefootnote{\fnsymbol{footnote}}
\setcounter{footnote}{1}

\input{Sections/Sec_Intro}
\input{Sections/Sec_BoGpOverview}
\input{Sections/Sec_ExistingCstrOptz}
\input{Sections/Sec_FrameworkCstrBo}
\input{Sections/Sec_BoVsQN}
\input{Sections/Sec_Conclusions}

\bmsection*{Acknowledgments}

The authors would like to recognize the Natural Sciences and Engineering Research Council of Canada and the Ontario Graduate Scholarship Program for the financial support they provided. The authors are also thankful for the helpful feedback provided by professors Masayuki Yano and Prasanth Nair from the University of Toronto while this work was being undertaken. 

% \bmsection*{Data availability statement}

\bmsection*{Conflict of interest}

The authors declare no potential conflict of interests.

\bmsection*{ORCID}

Andr\'e Marchildon: \url{https://orcid.org/0000-0001-6407-3987}

\bibliography{MyLibrary.bib}

\appendix
\input{Sections/Appendix_ExaLagMul}
\input{Sections/Appendix_TestCaseSol}

%\nocite{*}% Show all bib entries - both cited and uncited; comment this line to view only cited bib entries;

%\bmsection*{Author Biography}

%\begin{biography}{\includegraphics[width=76pt,height=76pt,draft]{empty}}{
%{\textbf{Author Name.} Please check with the journal's author guidelines whether
%author biographies are required. They are usually only included for
%review-type articles, and typically require photos and brief
%biographies for each author.}}
%\end{biography}

\end{document}

%% file: NewCommands.tex
\usepackage{bm}
\usepackage{bbm}
\usepackage[version=4]{mhchem}
\usepackage{siunitx}
\usepackage{longtable,tabularx}
\usepackage{xspace}

\usepackage{booktabs}
\usepackage{subcaption} 

\usepackage{graphicx}

\usepackage{float}
\floatstyle{plaintop}
\restylefloat{table}

\usepackage{tikz}

\usepackage{subcaption}

\usepackage{cancel}

\usepackage{algorithm,algpseudocode}

\usepackage{dsfont}

\usepackage{amsthm}

\usepackage{stackengine} % Required for the commanmd barbelow

\usepackage{lipsum}

\newcommand{\Table}[1]{\mbox{Table~\ref{#1}}}

\newcommand{\Fig}[1]{\mbox{Fig.~\ref{#1}}}
\newcommand{\Figs}[2]{\mbox{Figs.~\ref{#1}} and \ref{#2}}

\newcommand{\Sec}[1]{\mbox{Section~\ref{#1}}}

\newcommand{\Eq}[1]{\mbox{Eq.~(\ref{#1})}}
\newcommand{\Eqs}[2]{\mbox{Eqs.~(\ref{#1})} and (\ref{#2})}
\newcommand{\Eqss}[3]{\mbox{Eqs.~(\ref{#1})}, (\ref{#2}), and (\ref{#3})}

\newcommand{\bsym}[1]{\boldsymbol{#1}}

\newcommand{\mat}[1]{\ensuremath{\mathsf{#1}}}

\newcommand{\Idty}[0]{\mat{I}} % identity

\definecolor{Orange}{RGB}{255,69,0}
\definecolor{Red}{RGB}{235,40,25}
\definecolor{Green}{RGB}{34,139,34}
\newcommand\p[2]{\frac{\partial #1}{\partial #2}}
\newcommand\pp[3]{\frac{\partial^2 #1}{\partial #2 \partial #3}}

\newcommand{\order}[1]{\mathcal{O}\left(#1\right)}

\DeclareMathOperator{\mydiag}{diag}

\newcommand\lb[1]{\ensuremath{\underline{#1}}}
\newcommand\ub[1]{\ensuremath{\overline{#1}}}

\newcommand{\argmin}[0]{\ensuremath{\operatornamewithlimits{argmin}}}
\newcommand{\argmax}[0]{\ensuremath{\operatornamewithlimits{argmax}}}

\newcommand{\eigmin}[0]{\ensuremath{\lambda_{\text{min}}}}

\newcommand{\condmax}[0]{\kappa_{\max}}

\newcommand{\eg}[0]{e.g.\ }
\newcommand{\ie}[0]{i.e.\ }

\usepackage{amssymb}% http://ctan.org/pkg/amssymb
\usepackage{pifont}% http://ctan.org/pkg/pifont
\newcommand{\ignore}[1]{} % comment out large sections of code

\newcommand\yesnumber{\addtocounter{equation}{1}\tag{\theequation}}

\newcommand{\RNum}[1]{\uppercase\expandafter{\romannumeral #1\relax}}

\newcommand{\zero}[0]{\ensuremath{\textbf{0}}}

\newcommand{\one}[0]{\ensuremath{\textbf{1}}}

\newcommand{\bvec}[0]{\ensuremath{\bsym{b}}}

\newcommand{\fvec}[0]{\ensuremath{\bsym{f}}}
\newcommand{\gvec}[0]{\ensuremath{\bsym{g}}}
\newcommand{\hvec}[0]{\ensuremath{\bsym{h}}}

\newcommand{\kvec}[0]{\ensuremath{\bsym{k}}}

\newcommand{\rvec}[0]{\ensuremath{\bsym{r}}}

\newcommand{\uvec}[0]{\ensuremath{\bsym{u}}}

\newcommand{\xvec}[0]{\ensuremath{\bsym{x}}}
\newcommand{\yvec}[0]{\ensuremath{\bsym{y}}}

\newcommand{\betavec}[0]{\ensuremath{\bsym{\beta}}}

\newcommand{\muvec}[0]{\ensuremath{\bsym{\mu}}}
\newcommand{\sigmavec}[0]{\ensuremath{\bsym{\sigma}}}

\newcommand{\psivec}[0]{\ensuremath{\bsym{\psi}}}

\newcommand{\gammavec}[0]{\bsym{\gamma}}

\newcommand{\A}[0]{\mat{A}}

\newcommand{\G}[0]{\mat{G}}

\newcommand{\K}[0]{\mat{K}}
\newcommand{\M}[0]{\mat{M}}

\newcommand{\W}[0]{\mat{W}}
\newcommand{\X}[0]{\mat{X}}

\newcommand{\rdot}[0]{\ensuremath{\dot{r}}}

\newcommand{\rvecdot}[0]{\ensuremath{\dot{\rvec}}}

\newcommand{\fgrad}{\ensuremath{\bsym{f}_{\nabla}}}
\newcommand{\mgrad}{\ensuremath{\bsym{m}_{\nabla}}}

\newcommand{\onemod}[0]{\ensuremath{\check{\one}}}

\newcommand{\kG}[0]{\ensuremath{k_{\text{G}}}}
\newcommand{\kM}[0]{\ensuremath{k_{\text{M}{\frac{5}{2}}}}}

\newcommand{\kvecgrad}{\ensuremath{\bsym{k}_\nabla}}

\newcommand{\Kg}[0]{\ensuremath{\mat{K}_{\nabla}}}

\newcommand{\Kgdot}[0]{\ensuremath{\dot{\mat{K}}_{\nabla}}}

\newcommand{\etaKg}{\ensuremath{\eta_{\Kg}}}
\newcommand{\etaKgdot}{\ensuremath{\eta}_{\Kgdot}}

\newcommand{\Sigmag}[0]{\ensuremath{\Sigma_{\nabla}}}

\newcommand{\Pmat}[0]{\mat{P}}
\newcommand{\Pinv}[0]{\mat{P}^{-1}}

\newcommand{\sigK}{\ensuremath{\hat{\sigma}_{\K}}}
\newcommand{\sigKf}{\ensuremath{\hat{\sigma}_{\K,f}}}
\newcommand\muGP{\ensuremath{\tilde{\mu}}}
\newcommand\sigmaGP{\ensuremath{\tilde{\sigma}}}

\newcommand\muGPgvec[0]{\ensuremath{\tilde{\muvec}_g}}
\newcommand\sigmaGPgvec[0]{\ensuremath{\tilde{\sigmavec}_g}}

\newcommand\muGPhvec[0]{\ensuremath{\tilde{\muvec}_h}}
\newcommand\sigmaGPhvec[0]{\ensuremath{\tilde{\sigmavec}_h}}

\newcommand\Xdata[0]{\ensuremath{\X_{\text{data}}}}

\newcommand\Jdata[0]{\ensuremath{\bsym{J}_{\text{data}}}}

\newcommand\noisy[1]{\ensuremath{\widetilde{#1}}} % Noisy data

\newcommand\pdf[0]{\ensuremath{\theta_\text{pdf}}}
\newcommand\cdf[0]{\ensuremath{\theta_\text{cdf}}}

\newcommand\nlincong[0]{\ensuremath{n_{g,\text{lin}}}}
\newcommand\nlinconh[0]{\ensuremath{n_{h,\text{lin}}}}

\newcommand\nnlcg[0]{\ensuremath{n_{g,\text{nlin}}}}
\newcommand\nnlch[0]{\ensuremath{n_{h,\text{nlin}}}}

\newcommand\Lag[0]{\ensuremath{\mathcal{L}}}
\newcommand\augLag[0]{\ensuremath{\Lag_{\text{a}}}}
\newcommand\exaaugLag[0]{\ensuremath{\Lag_{\text{ea}}}}

\newcommand\xbest[0]{\ensuremath{\xvec_{\text{best}}}}

\newcommand\gtrc[0]{\ensuremath{g_{\text{trc}}}}
\newcommand\ubtrc[1]{\ensuremath{\ub{g}_{\text{trc}}^{#1}}}

\newcommand\gtrsig[0]{\ensuremath{g_{\text{tr}\sigmaGP_f}}}
\newcommand\ubtrsig[1]{\ensuremath{\ub{g}_{\text{tr}\sigmaGP_f}^{#1}}}

\newcommand\qEI[0]{\ensuremath{q_{\text{EI}}}}
\newcommand\qCEI[0]{\ensuremath{q_{\text{cEI}}}}

\newcommand\qUC[0]{\ensuremath{q_{\text{UC}}}}
\newcommand\qCUC[0]{\ensuremath{q_{\text{cUC}}}}

\newcommand\qExp[0]{\ensuremath{q_{\text{exp}}}}
\newcommand\qLag[0]{\ensuremath{q_{\Lag}}}

\newcommand\fbest[0]{\ensuremath{f_{\text{best}}}}

\newcommand\gvecAll[0]{\ensuremath{\gvec_{\text{all}}}}
\newcommand\hvecAll[0]{\ensuremath{\hvec_{\text{all}}}}

\newcommand\cZero[0]{\ensuremath{\textcolor{Red}{0}}}

\newcommand\cFive[0]{\ensuremath{\textcolor{Green}{5}}}
\newcommand\cTentyFive[0]{\ensuremath{\textcolor{Green}{25}}}

\usepackage{soul}

\definecolor{verylightblue}{rgb}{0.7,0.7,0.85}

%% file: Sections/Sec_Intro.tex
%!TEX root = ../Marchildon_ZinggCstrBayeOptz.tex

% ------------------------------------------
% New subsection
% ------------------------------------------
\section{Introduction} \label{Sec_Intro}

Bayesian optimizers are a popular and versatile class of numerical optimizers that have been used for various applications that include structures and aerodynamics \cite{mortished_aircraft_2016, march_gradient-based_2011,jim_bayesian_2021}. Bayesian optimizers use a probabilistic surrogate, which is commonly a Gaussian process (GP), to approximate a function of interest \cite{rasmussen_gaussian_2006}. The use of a probabilistic surrogate enables Bayesian optimizers to naturally use noisy function evaluations and to quantify the uncertainties in their predictions \cite{ameli_noise_2022}. Although Bayesian optimizers can also leverage gradients to minimize functions of interest with fewer function evaluations \cite{morris_bayesian_1993,ulaganathan_performance_2016}, they are commonly applied to unconstrained global optimization problems where gradients are not available \cite{shahriari_taking_2016,zhan_expected_2020}. Bayesian optimization has also had limited development for local optimization \cite{eriksson_scalable_2019} and problems with nonlinear constraints\cite{gardner_bayesian_2014,letham_constrained_2019}, which are the focus of this paper.

The use of gradients for a Bayesian optimizer enables its surrogate to more accurately approximate a function of interest \cite{han_improving_2013}. This is particularly useful when function evaluations are computationally expensive and when the gradients can be evaluated efficiently, such as with the adjoint method \cite{pironneau_optimum_1974,jameson_aerodynamic_1988}. However, the use of gradients generally makes the ill-conditioning problem of the covariance matrix significantly worse \cite{dalbey_efficient_2013,marchildon_non-intrusive_2023}. Fortunately, this problem can be efficiently addressed by preconditioning the covariance matrix and using a modest nugget \cite{marchildon_solution_2024}. This method guarantees that the Cholesky decomposition can be performed on a matrix with a condition number no greater than a user set threshold. 

Bayesian optimizers are effective global optimizers and have traditionally been applied to this class of problems since their surrogates can effectively approximate multimodal functions \cite{jones_efficient_1998,shahriari_taking_2016}. Furthermore, the probabilistic surrogate of the Bayesian optimizer enables exploration and exploitation tradeoffs in the multimodal parameter space \cite{nocedal_numerical_2006,zhan_expected_2020}. For local optimization, quasi-Newton optimizers are popular since they are robust and efficient \cite{byrd_global_1987}. However, it was recently demonstrated that a gradient-enhanced Bayesian optimizer can reach the same level of convergence as quasi-Newton optimizers for unimodal unconstrained problems and often do so with fewer function evaluations \cite{marchildon_efficient_2025}. The local optimization framework for the Bayesian optimizer uses a subset of the function and gradient evaluations to form the surrogate and to leverage its probabilistic surrogate to form a trust region. This local optimization framework is extended in this paper to problems with nonlinear constraints.

Bayesian optimizers have previously been applied to problems with nonlinear constraints, specifically nonlinear inequality constraints \cite{gardner_bayesian_2014,pourmohamad_bayesian_2021}. A popular method of handling nonlinear inequality constraints looks at the probability that each of the nonlinear constraints is satisfied \cite{gardner_bayesian_2014,letham_constrained_2019}. A separate GP is used to model each of the nonlinear constraints and the objective function. The acquisition function is then formed as the product of the expected improvement function and the probability that each of the nonlinear constraints is satisfied \cite{zhan_expected_2020}. However, the probability that a nonlinear equality constraint is satisfied is generally of measure zero, which limits this method to problems that only contain nonlinear inequality constraints. For many practical engineering optimization problems there are nonlinear equality constraints, such as a lift constraint when designing an aircraft \cite{chau_aerodynamic_2023,marchildon_gradient-enhanced_2024}.

Two similar methods were developed that enable Bayesian optimizers to handle both nonlinear inequality and equality constraints with the use of a non-exact augmented Lagrangian \cite{gramacy_modeling_2016,picheny_bayesian_2016}. A drawback of the first method is that it requires the use of a Monte Carlo method to evaluate the acquisition function. While the second method addressed this drawback, it was only tested on one two-dimensional optimization problem \cite{picheny_bayesian_2016}.

The goal of this paper is to enable Bayesian optimizers to solve nonlinearly-constrained optimization problems to a given tolerance using fewer function evaluations than quasi-Newton optimizers from SciPy and MATLAB. Two alternative methods for the Bayesian optimizer are presented and their performance is compared to quasi-Newton optimizers for several test cases containing two to 30 design variables.

The notation used in this paper is introduced in \Sec{Sec_Notation}. An overview of unconstrained Bayesian optimization is provided in \Sec{Sec_BoOverview} while constrained optimization is covered in \Sec{Sec_ExistingCstrOptz}. The two new methods that enable a Bayesian optimizer to handle nonlinear constraints are presented in \Sec{Sec_FrameworkCstrBo}. Test cases are introduced in \Sec{Sec_ConOptz_CstrTestCases} and used in \Sec{Sec_ConOptz_Studies} to compare the performance of the Bayesian optimizer with different methods of handling nonlinear constraints. The Bayesian optimizer with the two new methods is then compared to quasi-Newton optimizers in \Sec{Sec_BoVsQN}. The conclusions of this paper can be found in \Sec{Sec_Conclusions}. 

% ------------------------------------------
% New subsection
% ------------------------------------------
\section{Notation} \label{Sec_Notation}

Scalars are denoted by nonbolded lowercase sans-serif Latin and Greek letters while vectors are similarly denoted but bolded. For matrices, uppercase sans-serif Latin and Greek letters are used. The $ij$-th entry of the matrix $\X$ is $x_{ij}$, while its $i$-th column and $j$-th row are denoted by the vectors $\xvec_{:i}$ and $\xvec_{j:}$, respectively. All of the entries in the vectors $\zero$ and $\one$ are equal to zero and one, respectively. Integer quantities, such as the number of evaluation points and dimensions, are denoted by $n_x$ and $n_d$, respectively. When the nabla symbol $\nabla$ is used in the subscript of a matrix or vector, it indicates that it contains both function and gradient information. For example, $\fgrad$ is a vector containing all of the function and gradient evaluations of the function $f(\cdot)$. Finally, when a matrix is provided to the operator $\mydiag(\cdot)$, it outputs its diagonal entries as a vector, and when it is provided with a vector, it returns a diagonal matrix.

%% file: Sections/Sec_BoGpOverview.tex
%!TEX root = ../Marchildon_ZinggCstrBayeOptz.tex

% ------------------------------------------
% New subsection
% ------------------------------------------
\section{Overview of unconstrained Bayesian optimization} \label{Sec_BoOverview}

Bayesian optimizers require two ingredients: a probabilistic surrogate and an acquisition function. The probabilistic surrogate allows a function to be approximated and its uncertainty to be quantified. The acquisition function is constructed using the probabilistic surrogate and it is minimized to determine the next point in the parameter space where the function of interest should be evaluated. 

% ------------------------------------------
\subsection{Gradient-enhanced Gaussian Process} \label{Sec_BoOverview_GP}

The probabilistic surrogate is commonly a GP, which itself requires two components: a mean function, which is often simply taken to be a constant, and covariance function, which is commonly called a kernel function. Popular kernel functions include the Gaussian and Mat\'ern $\frac{5}{2}$ kernels:
\begin{alignat}{3}
	\kG(\xvec, \yvec; \gammavec) 
	&= \kG(\rvecdot)
	&&= e^{-\frac{1}{2} \| \rvecdot \|^2} \label{Eq_kern_Gaussian} \\
	\kM(\xvec, \yvec; \gammavec) 
	&= \kM(\rvecdot) 
	&&= \left(1 + \sqrt{3} \| \rvecdot \| + \| \rvecdot \|^2 \right) e^{- \sqrt{3} \| \rvecdot \|} \label{Eq_kern_Mat5f2},
%	\kRq(\xvec, \yvec; \gammavec, \alpha) 
%	&= \kRq(\rvecdot; \alpha) 
%	&&= \left(1 + \frac{ \| \rvecdot \|^2 }{2 \alpha} \right)^{-\alpha}, \label{Eq_kern_RatQd}
\end{alignat}
where the vector $\gammavec \in \mathbb{R}_+^{n_d}$ contains hyperparameters, and $\rvecdot$ is a dimensionless vector with its entries equal to $\rdot_i = \gamma_i (x_i - y_i) \, \forall \, i \in \{1, \ldots, n_d \}$. In this paper the Gaussian kernel is used, but other kernels, such as the Mat\'ern $\frac{5}{2}$ kernel, can also be used with the methods that are presented. 

To form the joint distribution and posterior for the gradient-enhanced GP we need the gradient-free and gradient-enhanced kernel matrices:
\begin{align} 
	\K(\X; \gammavec) 
	&=
	\begin{bmatrix}
		k(\xvec_{1:}, \xvec_{1:}; \gammavec) 	& k(\xvec_{1:}, \xvec_{2:}; \gammavec) 		& \ldots 	& k(\xvec_{1:}, \xvec_{n_x:}; \gammavec) \\
		k(\xvec_{2:}, \xvec_{1:}; \gammavec) 	& k(\xvec_{2:}, \xvec_{2:}; \gammavec) 		& \ldots 	& k(\xvec_{2:}, \xvec_{n_x:}; \gammavec) \\
		\vdots 	&	\vdots						& \ddots 	& \vdots \\
		k(\xvec_{n_x:}, \xvec_{1:}; \gammavec) 	& k(\xvec_{n_x:}, \xvec_{2:}; \gammavec) 	& \ldots 	& k(\xvec_{n_x:}, \xvec_{n_x:}; \gammavec)
	\end{bmatrix} \label{Eq_K} \\
	\Kg(\X; \gammavec) 
	&= 
	\begin{bmatrix}
		\K \left( \X; \gammavec \right) 				& \left(\p{\K}{y_1} \right)_\X			& \ldots 	& \left( \p{\K}{y_{n_d}} \right)_\X \\
		\left( \p{\K}{x_1} \right)_\X 	& \left( \pp{\K}{x_1}{y_1} \right)_\X 	& \ldots 	& \left( \pp{\K}{x_1}{y_{n_d}} \right)_\X \\ 
		\vdots 			& 	\vdots				& \ddots 	& \vdots \\
		\left( \p{\K}{x_{n_d}} \right)_\X 	& \left( \pp{\K}{x_{n_d}}{y_1} \right)_\X 	& \ldots 	& \left( \pp{\K}{x_{n_d}}{y_{n_d}} \right)_\X
	\end{bmatrix}, \label{Eq_Kg} 
\end{align} 
where $\left( \pp{\K}{x_i}{y_j} \right)_\X $ and the other submatrices of $\Kg$ are evaluated by applying the derivatives to each entry of the matrix $\K(\cdot; \gammavec)$ and then evaluating the resulting matrix at the rows of $\X$.

The joint distribution for the gradient-enhanced GP with noise-free function and gradient evaluations evaluated at the rows of $\X$, which are held in the vector $\fgrad(\X)$, along with the function evaluation at a point $\xvec'$ is 
%
% Page 16 of the Gaussian Processes for Machine Learning textbook
\begin{equation}
	\begin{bmatrix}
		\fgrad(\X) \\
		f(\xvec')
	\end{bmatrix}
	\sim \mathcal{N} \left(
	\begin{bmatrix}
		\zero \\
		0
	\end{bmatrix},
	\begin{bmatrix}
		\Sigmag & \sigK^2 \kvecgrad(\X; \xvec') \\
		\sigK^2 \kvecgrad(\xvec', \X) & \sigK^2 k(\xvec', \xvec')
	\end{bmatrix}
	\right),
\end{equation}
where the mean function is zero, $\Sigmag$ is the gradient-enhanced covariance matrix, and 
\begin{equation}
	\fgrad(\X) =
	\begin{bmatrix}
		\fvec(\X) \\
		\left(\p{\fvec}{x_1} \right)_{\X} \\
		\vdots \\
		\left(\p{\fvec}{x_{n_d}} \right)_{\X} \\
	\end{bmatrix}, \quad
	\kvecgrad(\X; \xvec') =
	\begin{bmatrix}
		\kvec(\X, \xvec') \\
		\left(\p{\kvec}{x_1} \right)_{(\X, \xvec')} \\
		\vdots \\
		\left(\p{\kvec}{x_{n_d}} \right)_{(\X, \xvec')}
	\end{bmatrix}.
\end{equation}
The gradient-enhanced covariance matrix is given by
\begin{equation}
	\Sigmag(\X; \sigK, \gammavec, \etaKg, \W) = \sigK^2 \left( \Kg(\X; \gammavec) + \etaKg \W \right),
\end{equation}
where $\etaKg \geq 0$, and $\W$ is a diagonal matrix with non-negative entries. The gradient-enhanced covariance matrix is commonly severely ill-conditioned. This problem can be addressed by using a preconditioning method \cite{marchildon_solution_2024} which ensures that the preconditioned matrix $\Pinv \Sigmag \Pinv$ has a condition number that is never greater than $\condmax > 1$, which is a user-defined parameter. When the function and gradient evaluations are noise-free, the preconditioning matrix $\Pmat$, the diagonal matrix $\W$, and the nugget $\etaKg$ are given by
\begin{align}
	\Pmat 
		&= \mydiag \left( \sqrt{ \mydiag \left( \Kg \right) } \right) \label{Eq_Pmat} \\
	\W 	&= \Pmat \Pmat \\
	\etaKgdot 
		&= \frac{ \max_{i} \sum_{j=1}^{n_d (n_d + 1)} \left| \Pmat^{-1} \Kg \Pmat^{-1} \right|_{ij} }{\condmax - 1}. \label{Eq_etaKgdot}
\end{align}
The preconditioning matrix must be adjusted when the function or gradient evaluations are noisy \cite{marchildon_solution_2024}. If the Gaussian kernel is used, then an upper bound for $\etaKg$ scales as $\order{n_x \sqrt{n_d}}$ \cite{marchildon_solution_2024}.

The mean and variance of the posterior are found by conditioning the prior of the GP on the noise-free observations, \ie the function and gradient evaluations at the rows of $\X$:
\begin{align}
	\muGP_f(\xvec') 
	&= m(\xvec') + \sigK^2 \, \kvecgrad(\xvec', \X) \Sigmag^{-1} \left(\noisy{\fgrad}(\X) - \mgrad(\X) \right) \label{Eq_muGp_w_grad} \\
	\sigmaGP_f^2(\xvec') 
	&= \sigK^2 \left( k(\xvec', \xvec') - \sigK^{-2} \, \kvecgrad(\xvec', \X) \Sigmag^{-1} \kvecgrad(\X, \xvec') \right). \label{Eq_sigmaGp_w_grad} 
\end{align}

% ------------------------------------------
\subsection{Marginal log-likelihood}

It is common to select the hyperparameters of the GP by finding the maximum of the likelihood function:
\begin{equation} \label{Eq_lkd}
	L(\gammavec, \beta, \sigK; \X, \fgrad, \etaKg, \W) 
	= \frac{e^{-\frac{ \left(\fgrad - \onemod \beta \right)^\top \Sigmag^{-1} \left( \fgrad - \onemod \beta \right) }{2}} }{ \left(2 \pi \right)^{\frac{n_x(n_d+1)}{2}} \sqrt{\det \left(\Sigmag \right) }},
\end{equation}
where $\onemod = [\one_{n_x}^\top, \zero_{n_x n_d}^\top]^\top$, and the hyperparameter $\beta$ is used as a constant mean function. There is a closed-form solution for the constant $\beta$ that maximizes the likelihood function. Furthermore, when all of the function and gradient evaluations are noise-free, there is also a closed form solution for $\sigK$:
\begin{align}
	\betavec
		&= \frac{ \onemod^\top \Sigmag^{-1} \fgrad}{ \onemod^\top \Sigmag^{-1} \onemod} \label{Eq_lkd_beta} \\
	\sigK^2 
		&= \frac{\left( \fgrad - \onemod \beta \right)^\top \left( \Kg + \etaKg \W \right)^{-1} \left( \fgrad - \onemod \beta \right)}{n_x (n_d+1)}. \yesnumber \label{Eq_lkd_sigK2_noisefree}
\end{align}
Taking the natural logarithm of the likelihood from \Eq{Eq_lkd}, substituting in $\sigK^2$ from \Eq{Eq_lkd_sigK2_noisefree}, and dropping the constant terms gives
\begin{equation} \label{Eq_ln_lkd_noisefree}
	\ln \left( L(\gammavec) \right) 
		= -\frac{n_x (n_d +1)}{2} \ln \left( \sigK^2 \right)
		- \frac{1}{2} \ln \left( \det \left( \Kg(\gammavec) + \etaKg \W \right) \right),
\end{equation}
which can be maximized with a numerical optimizer to select the hyperparameters $\gammavec$.

% ------------------------------------------
\subsection{Acquisition functions} \label{Sec_BoOverview_Acq}

Two popular acquisition functions for unconstrained Bayesian optimization are upper confidence (UC) and expected improvement (EI)\cite{jones_efficient_1998}:
\begin{align}
	q_{\text{UC}}(\xvec; \omega) 
		&= \muGP_f(\xvec) - \omega \sigmaGP_f(\xvec) \label{Eq_acq_UC} \\
	q_{\text{EI}}(\xvec; \fbest) 
		&= -\int_{-\infty}^{\fbest} (\fbest - f) \pdf \left(\frac{f - \muGP_f(\xvec)}{\sigmaGP_f(\xvec)} \right) df \nonumber \\
		&= -\left( \fbest - \muGP_f(\xvec) \right) \cdf \left(\frac{\fbest - \muGP_f(\xvec)}{\sigmaGP_f(\xvec)} \right) 
			- \sigmaGP_f(\xvec) \pdf \left(\frac{\fbest - \muGP_f(\xvec)}{\sigmaGP_f(\xvec)} \right), \label{Eq_acq_EI}
\end{align}
where a small value for $\omega > 0$ promotes exploitation while a large value promotes exploration, $\pdf$ and $\cdf$ are the probability and cumulative density functions, respectively, and $\fbest$ is the smallest objective function evaluation. Acquisition functions such as these are minimized to select the next point in the parameter space where the function of interest will be evaluated. Previously-used acquisition functions for nonlinearly-constrained Bayesian optimization are presented in \Sec{Sec_ConOptz_CurrentCstrBo}.

% ------------------------------------------
\subsection{Local optimization framework} \label{Sec_BoOverview_LocalOptzFramework}

This paper builds off of the local unconstrained Bayesian optimization framework from the same authors \cite{marchildon_efficient_2025} and a short overview of this framework is provided in this subsection. The parameters of this framework were selected by comparing numerical optimization results from different unimodal test cases with a range of dimensions.

First, the function and gradient evaluations from a subset of evaluation points are used to select the hyperparameters and to evaluate the posterior of the GP. This data region $\Xdata$ is centered around $\xbest$, \ie the evaluation point with the lowest merit function evaluation. The data region contains the 20 closest points to $\xbest$, in a Euclidean sense, including at least the last three evaluation points. The numerical results demonstrated that if too few or too many evaluation points are included in the data region, the Bayesian optimizer requires significantly more function evaluations to converge to the minimum.

The marginal log likelihood from \Eq{Eq_ln_lkd_noisefree} is maximized with a gradient-based optimizer to select the hyperparameters. Since the marginal log likelihood can be multimodal, the likelihood function is evaluated at 50 initial solutions and the optimization is started from the one with the highest likelihood evaluation. The preconditioning method is used to alleviate the ill-conditioning of the gradient-enhanced covariance matrix and the maximum condition number is set to $\condmax = 10^{10}$.

Two types of trust regions are used for the minimization of the acquisition function. The first, which is also commonly used by quasi-Newton optimizers, is simply a hypersphere around $\xbest$:
\begin{align}
	\gtrc(\xvec; \xbest) 
	&= \| \xvec - \xbest \|_2^2 \label{Eq_tr_circle_val} \\
	&\leq \ubtrc{j}, \nonumber
\end{align}
where $\ubtrc{j}$ is the upper bound for this trust region at the $j$-th optimization iteration. The second trust region leverages the uncertainty quantification from the posterior of the GP that is approximating the objective function:
\begin{align}
	\gtrsig(\xvec; \gammavec, \sigKf)
	&= \frac{\sigmaGP_f^2(\xvec; \gammavec, \sigKf)}{\sigKf^2} \label{Eq_tr_sigma_val} \\
	&= \left(k(\xvec,\xvec) - \sigKf^{-2} \, \kvecgrad(\X, \xvec)^\top \Sigmag^{-1} \kvecgrad(\X, \xvec) \right) \\
	&\leq \ubtrsig{j}, \nonumber
\end{align}
where $\ubtrsig{j}$ is the upper bound for this trust region at the $j$-th optimization iteration, and $k(\xvec,\xvec) = 1$ for the Gaussian and Mat\'ern $\frac{5}{2}$ kernels from \Eqs{Eq_kern_Gaussian}{Eq_kern_Mat5f2}, respectively. 
The trust region upper bounds $\ubtrc{j}$ and $\ubtrsig{j}$ are increased if progress is made and either of the trust regions was active at the minimized solution of the acquisition function. If progress is not made during two consecutive function evaluations, then the trust region upper bounds are decreased. Finally, if neither of these conditions are satisfied, the trust regions bounds are kept the same. The acquisition function is minimized with a quasi-Newton optimizer with both trust regions included as nonlinear constraints. Five independent optimizations are completed starting from different initial solutions and the final solution with the lowest acquisition function evaluation that satisfies the trust regions is kept.

%% file: Sections/Sec_ExistingCstrOptz.tex
%!TEX root = ../Marchildon_ZinggCstrBayeOptz.tex

% ------------------------------------------
% New section
% ------------------------------------------
\section{Overview of nonlinearly-constrained optimization}  \label{Sec_ExistingCstrOptz}

% ------------------------------------------
% New section
% ------------------------------------------
\subsection{Problem statement} \label{Sec_ConOptz_Intro}

A generalized constrained optimization problem is given by
\begin{align*}
	\min_{\lb{\xvec} \leq \xvec \leq \ub{\xvec}} f(\xvec) \quad \text{subject to} \quad \quad
	\A_g \xvec &\leq \bvec_g \yesnumber \label{Eq_optz_w_nlc} \\
	\A_h \xvec &= \bvec_h \\
	g_i(\xvec) & \leq 0 \quad \forall \, i \in \{1, \ldots, \nnlcg \} \\
	h_i(\xvec) & = 0 	\quad \forall \, i \in \{1, \ldots, \nnlch \},
\end{align*}
where $\A_g \in \mathbb{R}^{\nlincong \times n_d}$, $\bvec_g \in \mathbb{R}^{\nlincong}$, $\A_h \in \mathbb{R}^{\nlinconh \times n_d}$, and $\bvec_h \in \mathbb{R}^{\nlinconh}$ are for the linear inequality and equality constraints, respectively, while $g_i(\xvec) \, \forall \, i \in \{1, \ldots, \nnlcg \}$ and $h_i(\xvec) \, \forall \, i \in \{1, \ldots, \nnlch \}$ are the nonlinear inequality and equality constraints, respectively. Nonlinear inequality constraints given by $\lb{g} \leq g(\xvec) \leq \ub{g}$ can be recast as $\lb{g} - g(\xvec) \leq 0$ and $g(\xvec) - \ub{g} \leq 0$, and likewise for equality constraints. For the bound constraints we have $\lb{x}_i < \ub{x}_i \, \forall \, i \, \{1, \ldots, n_d \}$.

It is straightforward for an optimizer to check if a given point $\xvec$ in the parameter space satisfies the bound and linear constraints. Therefore, the objective function and nonlinear constraints, which may be expensive to evaluate, are generally only sampled at points $\xvec$ that satisfy the bound and linear constraints. In the following section, two methods that are commonly used by deterministic optimizers, such as quasi-Newton methods, are presented. These methods will then be adapted in \Sec{Sec_FrameworkCstrBo} to enable a Bayesian optimizer to handle nonlinear constraints. 

% ------------------------------------------
% New section
% ------------------------------------------
\subsection{Constrained optimization for deterministic optimizers} \label{Sec_ConOptz_MethodsExist}

% ------------------------------------------
% New subsection
% ------------------------------------------
\subsubsection{The $\ell_2$ penalty method} \label{Sec_ConOptz_MethodsExist_Penalty}

Since it is straightforward for most optimizers to ensure that the bound and linear constraints are satisfied exactly, they are omitted from this section. A merit function combines the objective and constraints into one single function. One such simple example is formed by adding an $\ell_2$ penalty for the nonlinear constraints to the objective function:
\begin{align} 
	J(\xvec) &= f(\xvec) + \rho \left( \| \hvec(\xvec) \|_2^2 + \| \gvec^+(\xvec) \|_2^2 \right), \label{Eq_merit_l2_penalty} 
\end{align}
where $\rho > 0$ and 
\begin{equation}
	g_i^+(\xvec) = \max(g_i(\xvec), 0) \quad \forall \, i \{1, \ldots, \nnlcg \}.
\end{equation}
This merit function is straightforward to implement and is at least twice differentiable if the objective function and constraints are as well. The $\ell_2$ penalty method is effective for penalizing large constraint violations, but it is ineffective for ensuring that the constraints are satisfied exactly since this requires that $\rho \rightarrow \infty$. 

% ------------------------------------------
% New subsection
% ------------------------------------------
\subsubsection{Augmented Lagrangian} \label{Sec_ConOptz_MethodsExist_AugLag}

To solve a constrained optimization problem without the ill-conditioning problems of the $\ell_2$ penalty method, a Lagrangian can be used instead:
\begin{equation} \label{Eq_Lag_base_nlc}
	L(\xvec; \psivec_h, \psivec_g) = f(\xvec) + \psivec_h^\top \hvec(\xvec) + \psivec_g^\top \gvec(\xvec),
\end{equation}
where $\psivec_h$ and $\psivec_g$ are the Lagrange multipliers for the equality and inequality constraints, respectively. The first-order optimality conditions, which are known as the KKT conditions \cite{nocedal_numerical_2006}, are given by
\begin{align}
	\p{L(\xvec)}{x_j} 
		= \p{f(\xvec)}{x_j} + \psivec_h^\top \p{\hvec(\xvec)}{x_j} + \psivec_g \p{\gvec(\xvec)}{x_j}
		&= 0 \quad \forall \, j \in \{1, \ldots, n_d \} \label{Eq_KKT_stationary} \\
	h_j(\xvec) 
		&= 0 \quad \forall \, j \in \{1, \ldots, \nnlch \} \label{Eq_KKT_h} \\
	g_j(\xvec)
		&\leq 0 \quad \forall \, j \in \{1, \ldots, \nnlcg \} \label{Eq_KKT_g_primal} \\
	\psivec_{g_j} 
		&\geq 0 \quad \forall \, j \in \{1, \ldots, \nnlcg \} \label{Eq_KKT_g_dual} \\
	\psivec_{g_j} g_j(\xvec) 
		&= 0 \quad \forall \, j \in \{1, \ldots, \nnlcg \} \label{Eq_KKT_comp_slack}.
\end{align}
The magnitudes of the Lagrange multipliers quantify how sensitive the constrained solution is to each of the constraints at an optimal point, \ie a solution to \Eq{Eq_optz_w_nlc}. A Lagrange multiplier is also needed for each of the bound and linear constraints that are active to ensure that the KKT conditions are satisfied.

The Lagrangian can be combined with the $\ell_2$ penalty method from \Sec{Sec_ConOptz_MethodsExist_Penalty} to form an augmented Lagrangian:
\begin{align} \label{Eq_merit_aug_Lag}
	\augLag(\xvec; \psivec_h, \psivec_g, \rho) 
	&= \left[ f(\xvec) + \psivec_h^\top \hvec(\xvec) + \psivec_g^\top \gvec(\xvec) \right] \\
	& \quad \quad + \rho \left[ \| \hvec(\xvec) \|_2^2 + \| \gvec(\xvec) \|_2^2 - \sum_{i = 1}^{\nnlcg} \min \left(0, \frac{\psi_{g_i}}{2 \rho} + g_i(\xvec) \right)^2 \right] \nonumber.
\end{align}
The $\ell_2$ penalty is efficient at penalizing portions of the parameter space that have a large infeasibility, but it is in general unable on its own to ensure that the nonlinear constraints are satisfied exactly. In contrast, the strength of the Lagrangian is to enable the nonlinear constraints to be satisfied exactly. It is straightforward to verify that \Eq{Eq_merit_aug_Lag} is equal to the objective evaluation if the KKT conditions are satisfied.

%Page 6 of https://epubs.siam.org/doi/pdf/10.1137/0806017
%Augmented Lagrangians combine the Lagrangian with, generally, an $\ell_2$ penalty. The Hestenes-Powell-Rockafellar augmented Lagrangian is given by
%
The Lagrange multipliers $\psivec_h$ and $\psivec_g$, along with the $\ell_2$ penalty parameter $\rho$, are selected prior to the minimization of the augmented Lagrangian. The Lagrange multipliers are then updated after each evaluation of the objective and nonlinear constraints. 
%The augmented Lagrangian is thus an $n_d$ dimensional problem with the entries of $\xvec$ as its only variables. 
Since the values of the Lagrange multipliers are generally different at each local minimum, having them fixed during the minimization of the augmented Lagrangian only enables a single minimum to be found. This limitation can be avoided by using an exact augmented Lagrangian, which is presented in the next subsection.

% ------------------------------------------
% New subsection
% ------------------------------------------
\subsubsection{Exact augmented Lagrangian} \label{Sec_ConOptz_MethodsExist_ExaAugLag}

A method that provides an unconstrained merit function whose minimum is a solution to the constrained optimization problem is referred to as an exact penalty method \cite{spedicato_exact_1994}. In this case, the solution to the constrained problem can be either a local or global minimum. The $\ell_2$ penalty from \Eq{Eq_merit_l2_penalty}, for example, requires that $\rho \rightarrow \infty$ in order to be an exact penalty method. In contrast, the augmented Lagrangian from \Eq{Eq_merit_aug_Lag} is only exact if the Lagrange multipliers satisfy the KKT conditions. However, the Lagrange multipliers are generally selected prior to the minimization of the augmented Lagrangian since solving for them simultaneously with $\xvec$ is more involved \cite{spedicato_exact_1994,nocedal_numerical_2006}. The values of the Lagrange multipliers that satisfy the KKT conditions are not known until a solution to the constrained optimization has been found. Consequently, the Lagrange multipliers are updated after each iteration and should converge to the values of the KKT conditions as the optimizer gets closer to a minimum \cite{nocedal_numerical_2006}. This makes the augmented Lagrangian from \Eq{Eq_merit_aug_Lag} an iterative penalty method rather than an exact penalty method.

Exact penalty functions can be either non-differentiable or differentiable \cite{dolgopolik_unified_2018}. Exact augmented Lagrangians, which are differentiable, were first developed to only handle equality constraints by Fletcher \cite{fletcher_class_1970}. Exact augmented Lagrangians were later generalized to also handle inequality constraints \cite{glad_multiplier_1979}. We consider the exact augmented Lagrangian presented by Di Pillo \cite{spedicato_exact_1994}, which uses the same augmented Lagrangian as \Eq{Eq_merit_aug_Lag}:
\begin{align} 
	\exaaugLag(\xvec; \rho) 
		&= \left[ f(\xvec) + \psivec_h^\top (\xvec) \hvec(\xvec) + \psivec_g^\top (\xvec) \gvec(\xvec) \right] \nonumber \\
		& \quad \quad + \rho \left[ \| \hvec(\xvec) \|_2^2 + \| \gvec(\xvec) \|_2^2 - \sum_{i = 1}^{\nnlcg} \min \left(0, \frac{\psi_{g_i}(\xvec) }{2 \rho} + g_i(\xvec) \right)^2 \right], \label{Eq_Lag_exa_aug}
\end{align}
where the Lagrange multipliers are now explicit functions of $\xvec$. To select the Lagrange multipliers the following function is minimized
\begin{equation} \label{Eq_exaLagMul}
	\Psi(\psivec_h, \psivec_g; \xvec) 
		= \| \nabla_{\xvec} \augLag(\xvec; \psivec_h, \psivec_g) \|^2 
		+ \alpha_1 \| \G(\xvec) \psivec_g \|^2 
		+ \alpha_2 w(\xvec) \left( \| \psivec_h \|_2^2 + \| \psivec_g \|_2^2 \right),
\end{equation}
where $\alpha_1 > 0$, $\alpha_2 > 0$, $\G(\xvec) = \mydiag(\gvec(\xvec))$, and
\begin{equation}
	w(\xvec) = \| \gvec^+(\xvec) \|_2^2 + \| h(\xvec) \|^2.
\end{equation}
The first two terms in $\Psi(\psivec_h, \psivec_g; \xvec)$ come from the KKT conditions introduced in \Sec{Sec_ConOptz_MethodsExist_AugLag}, while the last term adds a penalty when the constraints are not satisfied. Di Pillo \cite{spedicato_exact_1994} provides a closed-form solution for the Lagrange multipliers, which can be found in Appendix \ref{Sec_Appendix_ExaLagMul}.

% ------------------------------------------
% New subsection
% ------------------------------------------
\subsection{Existing methods for constrained Bayesian optimization} \label{Sec_ConOptz_CurrentCstrBo}

In general, the objective function and nonlinear constraints are all evaluated together. That is to say, whenever the objective function is evaluated, the nonlinear constraints are evaluated as well. It is important to minimize the number of function evaluations when these are computationally expensive. This motivates the use of a separate GP to approximate each of the nonlinear constraints instead of taking a less expensive but also less accurate approach, such as linearization, which is commonly used by quasi-Newton optimizers \cite{nocedal_numerical_2006}. For each of the GPs approximating a nonlinear constraint, along with the one approximating the objective function, their hyperparameters must be selected independently by maximizing the marginal log-likelihood from \Eq{Eq_ln_lkd_noisefree}. Once the hyperparameters have been selected, the following vectors can be formed, which hold the mean of the posteriors for the GPs approximating the inequality and equality constraints, respectively:
\begin{align}
	\muGPgvec(\xvec) 
		&= \left[ \muGP_{g_1}, \ldots, \muGP_{g_{\nnlcg}} \right]^\top \label{Eq_muGPgvec} \\
	\muGPhvec(\xvec) 
	&= \left[ \muGP_{h_1}, \ldots, \muGP_{h_{\nnlch}} \right]^\top. \label{Eq_muGPhvec}
\end{align}
Likewise, the vectors holding the variance of the GPs' posterior for the nonlinear constraints are given by
\begin{align}
	\sigmaGPgvec^2(\xvec) 
	&= \left[ \sigmaGP_{g_1}^2, \ldots, \sigmaGP_{g_{\nnlcg}}^2 \right]^\top \\
	\sigmaGPhvec^2(\xvec) 
	&= \left[ \sigmaGP_{h_1}^2, \ldots, \sigmaGP_{h_{\nnlch}}^2 \right]^\top.
\end{align}

The minimization of the acquisition function for a constrained optimization problem can be cast as
\begin{align*}
	\min_{\lb{\xvec} \leq \xvec \leq \ub{\xvec}} \, q(\xvec) \quad \text{subject to} \quad 
	\A_g \xvec &\leq \bvec_g \yesnumber \label{Eq_acq_min_w_lincon} \\
	\A_h \xvec &= \bvec_h,
\end{align*}
where the acquisition function is only dependent on $\xvec$ through the evaluations of the posterior of the GPs approximating the objective function and nonlinear constraints, \ie $q(\xvec) = q(\muGP_f(\xvec), \sigmaGP_f(\xvec), \muGPgvec(\xvec), \sigmaGPgvec(\xvec), \muGPhvec(\xvec), \sigmaGPhvec(\xvec))$. Trust regions can also be included in \Eq{Eq_acq_min_w_lincon} as nonlinear constraints.

A simple acquisition function for nonlinearly-constrained optimization problems involves using the $\ell_2$ penalty method presented in \Sec{Sec_ConOptz_MethodsExist_Penalty}:
\begin{equation} \label{Eq_acq_l2_pnlty}
	q_{\mu^2}(\xvec; \muGPhvec, \muGPgvec) = \| \muGPhvec(\xvec) \|_2^2 + \| \muGPgvec^+(\xvec) \|_2^2,
\end{equation}
where
\begin{equation}
	\muGP_{g_i}^+(\xvec) = \max \left( \muGP_{g_i}(\xvec), 0 \right) \quad \forall \, i \{1, \ldots, \nnlcg \}.
\end{equation}
This was one of the first methods used to apply a Bayesian optimizer to nonlinear constraints \cite{bjorkman_global_2000}. Barrier functions have also been used for Bayesian optimization but this only enables nonlinear inequality constraints to be considered \cite{pourmohamad_bayesian_2021}.

A popular method for Bayesian optimizers to handle nonlinear inequality constraints is to combine the expected improvement acquisition function $\qEI$ from \Eq{Eq_acq_EI} with the probability that each of the nonlinear constraints is satisfied \cite{schonlau_global_1998,gardner_bayesian_2014,durantin_analysis_2016}:
\begin{align}
	\qCEI(\xvec; \fbest) 
		&= \qEI(\xvec; \fbest) \prod_{i=1}^{\nnlcg} \mathbb{P}(g_i(\xvec) \leq 0) \nonumber \\
		&= \qEI(\xvec; \fbest) \prod_{i=1}^{\nnlcg} \cdf \left( 0; \muGP_{g_i}, \sigmaGP_{g_i} \right), \label{Eq_acq_cEI}
\end{align}
where $\mathbb{P}(\cdot)$ provides the probability that the condition it contains is satisfied, and $\fbest$ is the evaluation of the objective function with the smallest merit function. Only the inequality constraints can be considered since the probability that the equality constraints are satisfied is generally of measure zero. For noisy function evaluations, the probability that the constraints are satisfied can be adjusted to include a larger tolerance \cite{zhan_expected_2020}. % gelbart_bayesian_2014
Alternatively, the probability of feasibility can also be combined with the upper confidence acquisition function from \Eq{Eq_acq_UC} to get
\begin{align}
	\qCUC(\xvec; \omega) 
	&= \qUC(\xvec; \omega) \prod_{i=1}^{\nnlcg} \mathbb{P}(g_i(\xvec) \leq 0) \nonumber \\
	&= \qUC(\xvec; \omega) \prod_{i=1}^{\nnlcg} \cdf \left( 0; \muGP_{g_i}, \sigmaGP_{g_i} \right), \label{Eq_acq_cUC}
\end{align}
which will be compared to the use of $\qCEI(\xvec)$ and $q_{\mu^2}(\xvec)$ from \Eqs{Eq_acq_cEI}{Eq_acq_l2_pnlty}, respectively, in \Sec{Sec_ConOptz_Studies_ExistingMtds}.

%% file: Sections/Sec_FrameworkCstrBo.tex
%!TEX root = ../Marchildon_ZinggCstrBayeOptz.tex

% ------------------------------------------
% New subsection
% ------------------------------------------
\section{New constrained optimization framework} \label{Sec_FrameworkCstrBo}

% ------------------------------------------
% New subsection
% ------------------------------------------
\subsection{Probabilistic penalty to promote exploration} \label{Sec_ConOptz_MethodsNew_ProbPenalty}

The probability of feasibility from \Sec{Sec_ConOptz_CurrentCstrBo} cannot be applied to nonlinear equality constraints since they generally have a measure zero probability of being satisfied. This is the result of the following integration
\begin{align*}
	\mathbb{P}(h(\xvec) = 0) 
		&= \int_{\lb{h}}^{\ub{h}} \pdf(t; \muGP_h, \sigmaGP_h^2) dt \yesnumber \label{Eq_integral_pdf_h} \\
		&= \cdf(\ub{h}; \muGP_h, \sigmaGP_h) - \cdf(\lb{h}; \muGP_h, \sigmaGP_h) \\
		&= 0,
\end{align*}
where $\lb{h} = \ub{h} = 0$. The probability of feasibility for an equality constraint, which is based on the cumulative distribution function, is always measure zero for a probability distribution with a non-zero variance. However, this can be avoided by considering the probability density function instead. 

\begin{figure}[t!]
	\centering
	\begin{subfigure}[t]{0.45\textwidth}
		\includegraphics[width=\textwidth]{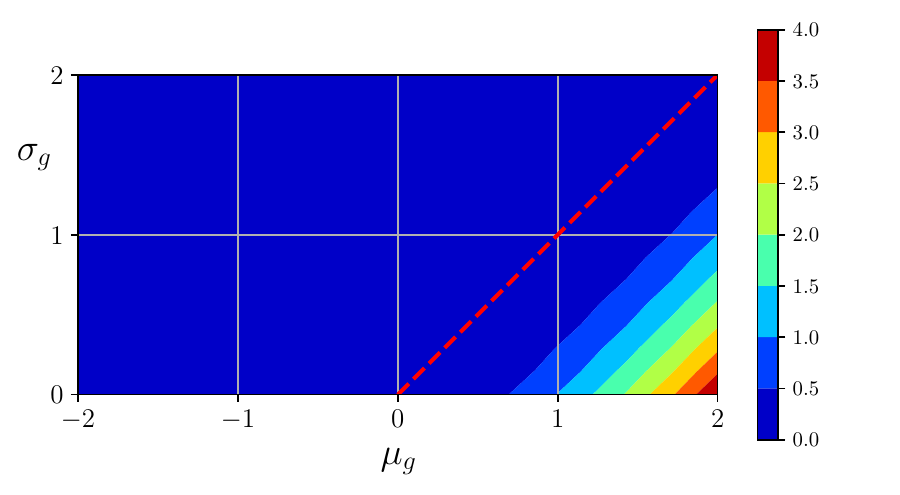}
		\caption{Penalty for nonlinear inequality constraints.}
		\label{Fig_AcqPenaltySigma_g}
	\end{subfigure}	
	\hspace{5pt}
	\begin{subfigure}[t]{0.45\textwidth}
		\includegraphics[width=\textwidth]{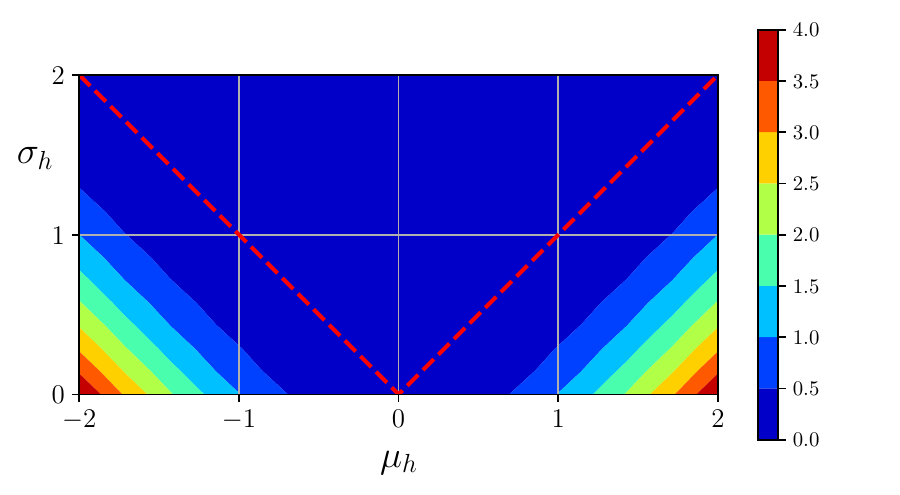}
		\caption{Penalty for nonlinear equality constraints.}
		\label{Fig_AcqPenaltySigma_h}
	\end{subfigure}	

	\caption{Acquisition penalty function $\qExp$ from \Eq{Eq_acq_qExp} for the nonlinear inequality and equality constraints. The red dashed line indicates $\mu = \sigma$, below this line $\qExp > 0$ while above it $\qExp = 0$.}
	\label{Fig_AcqPenaltySigma}
\end{figure}

The probability density function is always maximized with a mean and variance of zero, which results in the Dirac delta. However, the nonlinear constraints may not all be satisfied at the initial solution. In this case, it may not be possible to satisfy the nonlinear constraints at each iteration, \ie to have $\muGP_h = 0$, particularly if there are several constraints. We seek to find the variance $\sigmaGP_h^2$ that is the solution to the following maximization for a given mean $\muGP_h$:
\begin{equation} \label{Eq_pdf_optz}
	\sigmaGP_h^2 
		= \argmax_{\check{\sigma}_h^2} \pdf(t = 0; \muGP_h, \check{\sigma}_h^2),
\end{equation}
where $t=0$ since the equality constraint is satisfied for $h(\xvec) = 0$. 
%The probability density function is clearly maximized when $\muGP_h = 0$. 
To identify the value of $\sigmaGP_h^2$ that maximizes the probability density function, we start by calculating the following derivative:
\begin{align*}
	\p{\pdf(t; \muGP_h, \sigmaGP_h^2)}{\sigmaGP_h} 
		&= \p{}{\sigmaGP_h} \left( \frac{e^{-\frac{1}{2} \left( \frac{t - \muGP_h}{\sigmaGP_h} \right)^2 }}{\sigmaGP_h \sqrt{2 \pi}} \right) \\
		&= \left( \frac{(t - \muGP_h)^2}{\sigmaGP_h^4} -\frac{1}{\sigmaGP_h^2} \right) \frac{e^{-\frac{1}{2} \left( \frac{t - \muGP_h}{\sigmaGP_h} \right)^2 }}{\sqrt{2 \pi}} \\
		&= \left[ \left( \frac{t - \muGP_h}{\sigmaGP_h} \right)^2 - 1 \right] \frac{\pdf}{\sigmaGP_h}.
\end{align*}
For $t=0$ we have $\p{\pdf}{\sigmaGP_h} = 0$ when $\sigmaGP_h = | \muGP_h |$. To maximize the probability that an equality constraint is satisfied it is thus advantageous to have an acquisition function that penalizes large values of $|\muGP_h|$ and $|\muGP_h| - \sigmaGP_h$. The $\ell_2$ penalty function and the augmented Lagrangian already penalize $|\muGP_h|$. When the constraint is not satisfied, \ie $|\muGP_h| > 0$, penalizing large values of $|\muGP_h| - \sigmaGP_h$ thus promotes exploration to have a larger $\sigmaGP_h$, which increases the probability that the constraint is satisfied. This can be achieved with the following acquisition function
\begin{align}
	\qExp(\xvec)
		&= \sum_{i=1}^{\nnlcg} \max \left(\muGP_{g_i} - \sigmaGP_{g,_i} 0 \right)^2 
		+ \sum_{i=1}^{\nnlch} \left[ \max \left(\muGP_{h_i} - \sigmaGP_{h_i}, 0 \right)^2 + \max \left(-\muGP_{h_i} - \sigmaGP_{h_i}, 0 \right)^2 \right] \nonumber \\
		&= \sum_{i=1}^{\nnlcg} \max \left(\muGP_{g_i} - \sigmaGP_{g,_i} 0 \right)^2 
		+ \sum_{i=1}^{\nnlch} \max \left( \left| \muGP_{h_i} \right| - \sigmaGP_{h_i}, 0 \right)^2, \label{Eq_acq_qExp}
\end{align}
which adds a penalty when $\sigmaGP_h < |\muGP_h|$ for equality constraints and $\sigmaGP_g < \max(\muGP_g, 0)$ for inequality constraints. No penalty is added when $\sigmaGP_h \geq |\muGP_h|$, or $\sigmaGP_g \geq \max(\muGP_g, 0)$, since this would limit exploration. It is straightforward to verify that $\qExp(\xvec)$ is continuously differentiable. \Figs{Fig_AcqPenaltySigma_g}{Fig_AcqPenaltySigma_h} plot the acquisition function $\qExp(\xvec)$ when it is applied to a nonlinear inequality and equality constraint, respectively.

%--------------------------------------------------------------------------
%-------------------------------------------------------------------------- 

% ------------------------------------------
% New subsection
% ------------------------------------------
\subsection{Exact augmented Lagrangian} \label{Sec_ConOptz_MethodsNew_ExactAugLag}

An acquisition function can be formed using the augmented Lagrangian presented in \Sec{Sec_ConOptz_MethodsExist_ExaAugLag}. In order for the Lagrange multipliers to satisfy the KKT conditions, which are listed in \Sec{Sec_ConOptz_MethodsExist_AugLag}, the bound, linear, and nonlinear constraints must all be considered by the exact augmented Lagrangian. The following vectors are formed that contain all of the inequality and equality constraints:
\begin{align}
	\gvecAll(\xvec) 
		&= \left[ \left(\lb{\xvec} - \xvec \right)^\top, \left( \xvec - \ub{\xvec} \right)^\top, \left( \A_g \xvec -\bvec_g \right)^\top, \muGPgvec^\top(\xvec) \right]^\top \label{Eq_gvecAll} \\
	\hvecAll(\xvec) 
	&= \left[ \left( \A_h \xvec - \bvec_h \right)^\top, \muGPhvec^\top(\xvec) \right]^\top, \label{Eq_hvecAll}
\end{align}
where the mean of the posterior for the GPs $\muGPgvec$ and $\muGPhvec$ are from \Eqs{Eq_muGPgvec}{Eq_muGPhvec}, respectively. The lengths of the vectors $\gvecAll$ and $\hvecAll$ are $n_{g,\text{all}} = 2 n_d + \nlincong + \nnlcg$ and $n_{h,\text{all}} = \nlinconh + \nnlch$, respectively. Selecting the Lagrange multipliers $\psivec_h(\xvec)$ and $\psivec_g(\xvec)$ by minimizing $\Psi$ from \Eq{Eq_exaLagMul} requires inverting a symmetric positive definite matrix at a cost of $\order{(n_{g,\text{all}} + n_{h,\text{all}})^3}$, as shown in Appendix \ref{Sec_Appendix_ExaLagMul}. Since the Lagrange multipliers depend on $\xvec$, they must be updated each time the surrogate models are evaluated at a new point in the parameter space. This can be expensive when there are a large number of constraints, including linear and bound constraints. To reduce this computational cost, the vector $\gvecAll$ can be shortened by including only inequality constraints that are active or nearly active. Specifically, the vector $\gvecAll$ is formed with all of the inequality constraints when the acquisition function is sampled at a new point $\xvec$. Then, the entries in $\gvecAll$ that are smaller than $-0.1$ are removed from the vector. A small negative number is used to ensure that the nonlinear inequality constraints are included in the acquisition function whenever they are active and also for a small region where they are not active. This helps keep the acquisition function smooth in sections of the design space where the inequality constraints transition from being inactive to active.

The acquisition function formed with the exact augmented Lagrangian from \Sec{Sec_ConOptz_MethodsNew_ExactAugLag} is given by
\begin{align} 
	\qLag(\xvec; \rho) 
		&= \left[ \muGP_f(\xvec) + \psivec_h^\top (\xvec) \hvecAll(\xvec) + \psivec_g^\top (\xvec) \gvecAll(\xvec) \right] \nonumber \\
		& \quad \quad + \rho \left[ \| \hvecAll(\xvec) \|_2^2 + \| \gvecAll(\xvec) \|_2^2 - \sum_{i = 1}^{\nnlcg} \min \left(0, \frac{\psi_{g_i}(\xvec) }{2 \rho} + g_{i,\text{all}}(\xvec) \right)^2 \right]. \label{Eq_acq_Lag_exa_aug}
\end{align}
The Lagrange multipliers $\psivec_h$ and $\psivec_g$ are selected by minimizing $\Psi$ from \Eq{Eq_exaLagMul} with the coefficients $\alpha_1 = \alpha_2 = 100$. The vectors $\hvecAll$ and $\gvecAll$ from \Eqs{Eq_gvecAll}{Eq_hvecAll}, respectively, are used along with $\muGP_f(\xvec)$, \ie the mean of the posterior of the GP approximating the objective function. The steps needed to apply this method with a Bayesian optimizer to solve nonlinearly-constrained optimization problems are summarized in Algorithm \ref{Alg_ExaAugLag}.

%-------------------------------------------------------------------------- 
\begin{algorithm}[t!]
	\caption{Acquisition function minimization with an exact augmented Lagrangian}
	\label{Alg_ExaAugLag}
	\begin{algorithmic}[1]
		\Statex{\textbf{Required:} The matrix of evaluation points in the data region $\Xdata$ detailed in \Sec{Sec_BoOverview_LocalOptzFramework} and the merit function evaluation $\Jdata$ at these data points.}
		\Statex{\textbf{Select:} $0 < \epsilon_{g} (-10^{-1})$, $0 < \rho_1 (100)$, and $0 < \rho_2 (100)$, with default values indicated in parentheses}
		\State{From \Eq{Eq_gvecAll}: $\gvecAll(\xvec) = \left[ \left(\lb{\xvec} - \xvec \right)^\top, \left( \xvec - \ub{\xvec} \right)^\top, \left( \A_g \xvec -\bvec_g \right)^\top, \muGPgvec^\top(\xvec) \right]^\top$}
		\State{Remove all entries in $\gvecAll$ that do not satisfy $g_{i,\text{all}} \geq \epsilon_{g}$}
		\State{From \Eq{Eq_hvecAll}: $\hvecAll(\xvec) = \left[ \left( \A_h \xvec - \bvec_h \right)^\top, \muGPhvec^\top(\xvec) \right]^\top$}
		\State{Evaluate the Lagrangian multipliers with \Eq{Eq_ExaLagMul_sol} using $\gvecAll(\xvec)$, $\hvecAll(\xvec)$, and $\muGP_f(\xvec)$}
		\State{Evaluate the exact augmented Lagrangian acquisition function $\qLag(\xvec; \rho_1)$ with \Eq{Eq_acq_Lag_exa_aug}}
		\State{Form the complete acquisition function $q(\xvec) = \qLag(\xvec) + \rho_2 \qExp(\xvec)$, with $\qExp$ from \Eq{Eq_acq_qExp}}
		\State{Minimize the acquisition function $q(\xvec)$ with a gradient-based optimizer}
		\State{\textbf{Return:} The next point in the parameter space where the objective and nonlinear constraints will be evaluated: $\xvec_{\text{sol}}$.}
	\end{algorithmic}
\end{algorithm} 
%-------------------------------------------------------------------------- 

% ------------------------------------------
% New subsection
% ------------------------------------------
\subsection{Strong enforcement} \label{Sec_ConOptz_MethodsNew_StrongEnf}

All of the methods that have been presented thus far for nonlinear constrained optimization have included the nonlinear constraints in the acquisition function. This enables the minimization of the acquisition function to be formulated as \Eq{Eq_acq_min_w_lincon}, where the only constraints the acquisition function minimizer has to consider are the bound and linear constraints, and potentially the trust regions from \Sec{Sec_BoOverview_LocalOptzFramework}, which are nonlinear constraints. 

In this section the minimization of the acquisition function is reformulated to have additional constraints, which are the means of the posterior of the GPs that are approximating the nonlinear constraints:
\begin{align*}
	\xvec^{i+1} 
		= \argmin_{\lb{\xvec} \leq \xvec \leq \ub{\xvec}} \, q(\xvec) \quad \text{subject to} \quad \quad \quad
		\A_g \xvec &\leq \bvec_g \yesnumber \label{Eq_acq_min_nlc_strong_baseline} \\
		\A_h \xvec &= \bvec_h \\
		\gtrc(\xvec) &\leq \ubtrc{j} \\
		\gtrsig(\xvec) &\leq \ubtrsig{j} \\
		\muGP_{g_i}(\xvec) & \leq 0 \quad \forall \, i \in \{1, \ldots, \nnlcg \} \\
		\muGP_{h_i}(\xvec) & = 0 	\quad \forall \, i \in \{1, \ldots, \nnlch \},
\end{align*}
where $\gtrc(\xvec)$ and $\gtrsig(\xvec)$ are the circular and $\sigma$ trust regions from \Eqs{Eq_tr_circle_val}{Eq_tr_sigma_val}, respectively, while $\ubtrc{j}$ and $\ubtrsig{j}$ are their upper bounds for the $j$-th optimization iteration, respectively. 

There are two limitations to solving \Eq{Eq_acq_min_nlc_strong_baseline} at all optimization iterations. First, the mean of the posterior of the GP approximating the nonlinear constraints will not be accurate if there are too few function evaluations. Second, if the starting point is not feasible, it may not be possible to ensure that the nonlinear constraints $\muGP_{g_i}(\xvec) \leq 0 \, \forall \, i \in \{1, \ldots, \nnlcg \}$ and $\muGP_{h_i}(\xvec) = 0 \, \forall \, i \in \{1, \ldots, \nnlch \}$ can be satisfied, particularly if the trust regions are small. To address this, the strong enforcement of the nonlinear constraints is implemented in three stages. 

\vspace{\baselineskip}
\noindent\textbf{Stage 1}

In the first stage, \ie when $n_x < 10$, the constraints are not enforced explicitly. The $\ell_2$ penalty acquisition function $q_{\mu^2}(\xvec)$ from \Eq{Eq_acq_l2_pnlty} is used to handle the nonlinear constraints, which is efficient at reducing the infeasibility when it is large, \ie $\max(|\muGPhvec|) \gg 0$ or $\max(\muGPgvec) \gg 0$. The acquisition function minimization problem is thus
\begin{align*}
	\xvec^{i+1} = \argmin_{\xvec} q(\xvec) \quad \text{s.t.} \quad \quad 
	\gtrc(\xvec) 	&\leq \ubtrc{j} \yesnumber \label{Eq_acq_optz_uncon} \\
	\gtrsig(\xvec) 	&\leq \ubtrsig{j},
\end{align*}
where a suitable acquisition function would be $q(\xvec) = \qUC(\xvec) + 10^2 q_{\mu^2}(\xvec) + 10^2\qExp(\xvec)$, where $\qUC$, $q_{\mu^2}$, and $\qExp$ come from \Eqss{Eq_acq_UC}{Eq_acq_l2_pnlty}{Eq_acq_qExp}, respectively.

\vspace{\baselineskip}
\noindent\textbf{Stage 2}

\begin{figure}[t!]
	\centering
	\begin{subfigure}[t]{0.45\textwidth}
		\includegraphics[width=\textwidth]{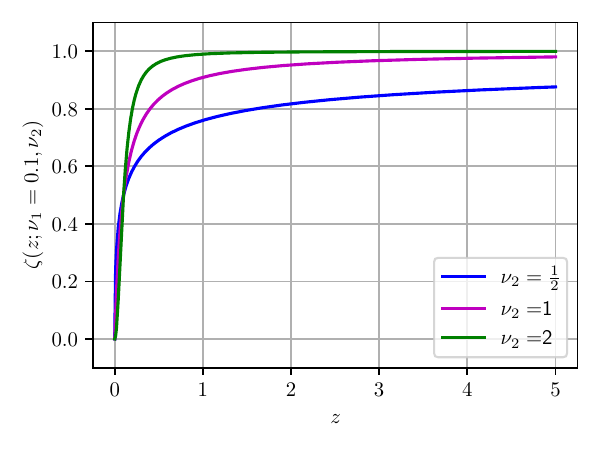}
		\caption{Linear scale.}
		\label{Fig_sigmoid_linear}
	\end{subfigure}	
	\hspace{5pt}
	\begin{subfigure}[t]{0.45\textwidth}
		\includegraphics[width=\textwidth]{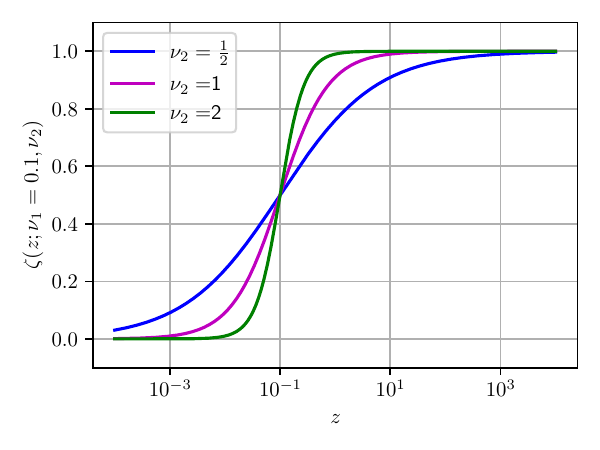}
		\caption{Logarithmic scale.}
		\label{Fig_sigmoid_log}
	\end{subfigure}	
	\caption[Plots for the sigmoid function that is used to reduce the upper bound for the nonlinear constraints when they are enforced using the strong method.]{Plots for the sigmoid function from \Eq{Eq_sigmoid_strong_enf} that is used to reduce the upper bound for the nonlinear constraints when they are enforced using the strong method detailed in Algorithm~\ref{Alg_StrongEnf}.}
	\label{Fig_sigmoid}
\end{figure}

The second stage starts when $n_x \geq 10$ and continues to be used while the feasibility remains larger than a set threshold. As a continuously differentiable measure of the feasibility we use $q_{\mu^2}(\xvec)$ from \Eq{Eq_acq_l2_pnlty}, which depends on $\muGPgvec(\xvec)$ and $\muGPhvec(\xvec)$. The mean of the posterior for the GPs approximating the nonlinear constraints are used instead of the nonlinear constraint evaluations, which are available at $\xbest$, since these could be noisy evaluations. For noise-free evaluations, there should be little difference between the nonlinear constraints and their approximations $\muGPgvec$ and $\muGPhvec$ at all of points in the parameter space where the constraints have been evaluated. During the second stage, a constraint is placed on $q_{\mu^2}(\xvec)$, while for stage 3 each of the nonlinear constraints are considered individually. By default, the same acquisition function $q(\xvec)$ as stage 1 is used. The formulation for the minimization of the acquisition function during the second stage of the strong enforcement of the nonlinear constraints is given by
\begin{align*}
	\xvec^{i+1} 
	= \argmin_{\lb{\xvec} \leq \xvec \leq \ub{\xvec}} \, q(\xvec) \quad \text{subject to} \quad \quad \quad
	\A_g \xvec	 	 &\leq \bvec_g \yesnumber \label{Eq_acq_min_nlc_strong_v2} \\
	\A_h \xvec 	 	 &= \bvec_h \\
	\gtrc(\xvec) 	 &\leq \ubtrc{j} \\
	\gtrsig(\xvec) 	 &\leq \ubtrsig{j} \\
	q_{\mu^2}(\xvec) & \leq \ub{q}_{\mu^2},
\end{align*}
where $\ub{q}_{\mu^2} > 0$ is the upper bound for $q_{\mu^2}(\xvec)$. It is selected with
\begin{equation} \label{Eq_ub_q_mu2}
	\ub{q}_{\mu^2}(\xbest) = \zeta \left(q_{\mu^2}(\xbest) ; \nu_1, \nu_2 \right) q_{\mu^2}(\xbest),
\end{equation}
where $\nu_1, \nu_2 > 0$ and 
\begin{equation} \label{Eq_sigmoid_strong_enf}
	\zeta(z; \nu_1, \nu_2) = \frac{(\nu_1 z)^{\nu_2}}{(\nu_1 z)^{\nu_2} + 1},
\end{equation}
which is a sigmoid function since $\zeta: \mathbb{R} \rightarrow [0,1]$ for $\nu_1, \nu_2 \in \mathbb{R}_+$ \footnote{It is in fact the logistic function: $f_{\text{logistic}}(x; \nu_2) = \frac{1}{1 + e^{-\nu_2 x}}$ with $x = \ln(\nu_1 z)$, which is type of sigmoid function.}. \Fig{Fig_sigmoid} plots the sigmoid function for $\nu_1 = 10$ and different values of $\nu_2$. The value of $z$ for $\zeta(z;\nu_1, \nu_2)$ where $\zeta(z) = \frac{1}{2}$ is given by $\nu_1^{-1}$.

%--------------------------------------------------------------------------
%-------------------------------------------------------------------------- 

\vspace{\baselineskip}
\noindent\textbf{Stage 3}

The third stage for the strong enforcement of the nonlinear constraints starts when the $\ell_2$ penalty is sufficiently small, \eg $q_{\mu^2} < 1$. In this case, each of the nonlinear constraints is considered independently and the following constrained minimization of the acquisition function is solved:
\begin{align*}
	\xvec^{i+1} 
	= \argmin_{\lb{\xvec} \leq \xvec \leq \ub{\xvec}} \, q(\xvec) \quad \text{subject to} \, \quad \quad \quad \quad \quad \quad
	\A_g \xvec &\leq \bvec_g \label{Eq_acq_min_nlc_strong_v3} \\
	\A_h \xvec &= \bvec_h \\
	\gtrc(\xvec) - \ubtrc{j} &\leq 0 \yesnumber \\
	\gtrsig(\xvec) - \ubtrsig{j} &\leq 0 \\
	\muGP_{g_i}(\xvec) & \leq \ub{\muGP}_{g_i} \quad \forall \, i \in \{1, \ldots, \nnlcg \} \\
	-\ub{\muGP}_{h_i} \leq \muGP_{h_i}(\xvec) & \leq \ub{\muGP}_{h_i}	\quad \forall \, i \in \{1, \ldots, \nnlch \},
\end{align*}
where the acquisition function is by default the same one used in the prior stages and
\begin{alignat}{2}
	\ub{\muGP}_{g_i}(\xbest) 
		&= \zeta \left( \muGP_{g_i}^+(\xbest), 0; \nu_1, \nu_2 \right) \muGP_{g_i}^+(\xbest) \quad 
		&&\forall \, i \{1, \ldots, \nnlcg \} \label{Eq_ub_muGP_g} \\
	\ub{\muGP}_{h_i}(\xbest) 
		&= \zeta \left( \muGP_{h_i}(\xbest); \nu_1, \nu_2 \right) \muGP_{h_i}(\xbest) \quad 
		&&\forall \, i \{1, \ldots, \nnlch \} \label{Eq_ub_muGP_h}.
\end{alignat}
The three stages for the strong enforcement of the nonlinear constraints are summarized in Algorithm~\ref{Alg_StrongEnf}. 

%-------------------------------------------------------------------------- 
\begin{algorithm}[t!]
	\caption{Acquisition function minimization with strong nonlinear constraint enforcement}
	\label{Alg_StrongEnf}
	\begin{algorithmic}[1]
		\Statex{\textbf{Required:} The matrix of evaluation points in the data region $\Xdata$ detailed in \Sec{Sec_BoOverview_LocalOptzFramework}, the merit function evaluation $\Jdata$ at these data points, and the mean and variances of the posterior for the GPs approximating the objective function and nonlinear constraints.}
		\Statex{\textbf{Select:} $0 < \epsilon_{\ell_2} (1)$, $0 < \nu_1 (10)$, and $0 < \nu_2 (1)$, $q (\qUC + 10^2 q_{\mu^2} + 10^2\qExp)$, with default values indicated in parentheses.}
		\If{$n_x < 10$}
		\Comment{Stage 1}
		\State{Solve the constrained acquisition function from \Eq{Eq_acq_optz_uncon}}
		\Else
		\State{Identify $\xbest$, \ie $\xbest = \left( \Xdata \right)_{i^*:}$, where $i^* = \argmin_i \left(\Jdata \right)_i$}
		\State{From \Eq{Eq_acq_l2_pnlty}: $q_{\mu^2}(\xbest) = \| \muGPhvec(\xbest) \|_2^2 + \| \muGPgvec^+(\xbest) \|_2^2$}
		\If{$q_{\mu^2}(\xbest) \geq \epsilon_{\ell_2}$}
		\Comment{Stage 2}
		\State{From \Eq{Eq_ub_q_mu2}: $\ub{q}_{\mu^2} = \zeta \left(q_{\mu^2}(\xbest) ; \nu_1, \nu_2 \right) q_{\mu^2}(\xbest)$}
		\State{Solve the constrained acquisition function from \Eq{Eq_acq_min_nlc_strong_v2}}
		\Else 
		\Comment{Stage 3}
		\State{From \Eq{Eq_ub_muGP_g}: $\ub{\muGP}_{g_i}(\xbest) = \zeta \left( \muGP_{g_i}(\xbest); \nu_1, \nu_2 \right) \muGP_{g_i}(\xbest) \, \forall \, i \{1, \ldots, \nnlcg \}$}
		\State{From \Eq{Eq_ub_muGP_h}: $\ub{\muGP}_{h_i}(\xbest) = \zeta \left( \muGP_{h_i}(\xbest); \nu_1, \nu_2 \right) \muGP_{h_i}(\xbest) \, \forall \, i \{1, \ldots, \nnlch \}$}
		\State{Solve the constrained acquisition function from \Eq{Eq_acq_min_nlc_strong_v3}}
		\EndIf
		\EndIf
		\State{\textbf{Return:} $\xvec_{\text{sol}}$, which is the next point in the parameter space where the objective and nonlinear constraints will be evaluated.}
	\end{algorithmic}
\end{algorithm}
%-------------------------------------------------------------------------- 

% ------------------------------------------
% New section
% ------------------------------------------
\section{Constrained test cases} \label{Sec_ConOptz_CstrTestCases}

\begin{figure}[t!]
	\centering
	\begin{subfigure}[t]{0.328\textwidth}
		\includegraphics[width=\textwidth]{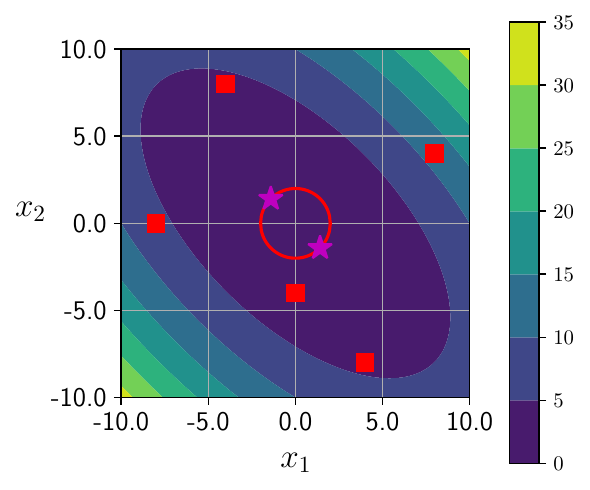}
		\caption{Constrained quadratic function from \Eq{Eq_min_quad_nlc_g}.}
		\label{Fig_BoUnconTestCases_quad_fun}
	\end{subfigure}	
	\begin{subfigure}[t]{0.328\textwidth}
		\includegraphics[width=\textwidth]{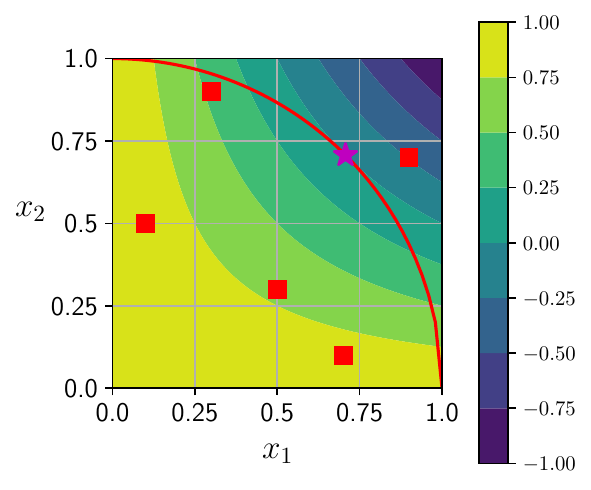}
		\caption{Constrained product function from \Eq{Eq_min_prod_nlc_h}.}
		\label{Fig_BoUnconTestCases_prod_fun}
	\end{subfigure}	
	\begin{subfigure}[t]{0.328\textwidth}
		\includegraphics[width=\textwidth]{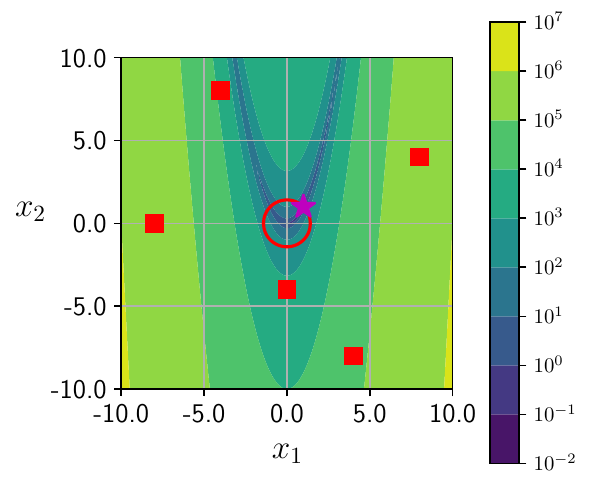}
		\caption{Constrained Rosenbrock function from \Eq{Eq_min_nlc_g_Rosen}.}
		\label{Fig_BoUnconTestCases_Rosen_fun}
	\end{subfigure}	
	\caption[Constrained test cases with a Latin hypercube sampling used to select the starting points of the optimizer.]{Plots for the constrained test cases with the red squares indicating the starting points for the optimizers, which were selected with a Latin hypercube sampling, and the magenta star indicating the solutions for the constrained minimization test cases. The red line denotes the constraint.}
	\label{Fig_BoUnconTestCases}
\end{figure}

Three constrained test cases are presented in this section. These will be used in \Sec{Sec_ConOptz_Studies} to test the Bayesian optimizer with different acquisition functions and settings for nonlinear constrained optimization. Then, in \Sec{Sec_BoVsQN} the constrained Bayesian optimizer will be compared to quasi-Newton optimizers with the test cases from this section. The following objective functions are used along with a nonlinear constraint to benchmark the optimizers:
\begin{alignat}{3} 
	\text{Centered quadratic function:} \quad
		& f_{\text{quad}}(\xvec; \A) 
		&&= \xvec^\top \A \xvec \label{Eq_Quadratic_centered_fun} \\
	\text{Product function:}  \quad
		& f_{\text{prod}}(\xvec) 
		&&= 1 - n_d^{n_d / 2} \prod_{i=1}^{n_d} x_i \label{Eq_prod_fun} \\
	\text{Rosenbrock function:} \quad
		&f_{\text{Rosen}}(\xvec)  
		&&= \sum_{i=1}^{n_d-1} \left[ a \left(x_{i+1} - x_i^2 \right)^2 + \left( 1 - x_i \right)^2 \right], \label{Eq_Rosenbrock}
\end{alignat}
where $a>0$ for the Rosenbrock function and $\A$ is a matrix with $a_{ij} = \frac{1}{10} e^{-\frac{1}{2}(i-j)^2} \, \forall i,j \in \{1, \ldots, n_d \}$. The three nonlinearly-constrained optimization problems are 
\begin{alignat}{5} 
	\min_{-10 \leq \xvec \leq 10} \,
	&f_{\text{quad}}(\xvec; \A) -4 \eigmin(\A) \quad 
	&&\text{subject to} \quad \| \xvec \|_2^2
	&&&&\geq 4 \label{Eq_min_quad_nlc_g} \\
	\min_{\phantom{-1}0 \leq \xvec \leq 1 \phantom{0}} \,
	&f_{\text{prod}}(\xvec) \quad 
	&&\text{subject to} \quad \| \xvec \|_2^2
	&&&&= 1 \label{Eq_min_prod_nlc_h} \\
	\min_{-10 \leq \xvec \leq 10} \,
	&f_{\text{Rosen}}(\xvec) \quad 
	&&\text{subject to} \quad \| \xvec \|_2^2
	&&&&\leq n_d. \label{Eq_min_nlc_g_Rosen}
\end{alignat}
The objective evaluates to zero for the solution to these three minimization problems. For \Eq{Eq_min_quad_nlc_g} there are two equivalent solutions at $\xvec = \pm 4 \uvec_{\min}$, where $\uvec$ is the unit eigenvector of $\A$. These solutions are derived in Appendix \ref{Sec_Appendix_TestCase_SolNlcQuad}. For \Eq{Eq_min_prod_nlc_h} its unique solution is $\xvec = n_d^{-1/2} \one$, as demonstrated in Appendix \ref{Sec_Appendix_TestCase_SolNlcProd}. The nonlinear inequality constraint of \Eq{Eq_min_nlc_g_Rosen} is active at its solution, \ie $\| \xvec \|_2^2 = n_d$, but it has the same solution as the unconstrained problem. The constraint for \Eq{Eq_min_nlc_g_Rosen} does nonetheless impact the minimization problem since most starting points for the optimizer do not satisfy the constraint. The three constrained test cases can be seen in \Fig{Fig_BoUnconTestCases} for the $n_d = 2$ case.

%--------------------------------------------------------------------------
%--------------------------------------------------------------------------

To compare the progress of the optimizers the exact augmented Lagrangian from \Eq{Eq_Lag_exa_aug} is used as the merit function with $\rho = 100$, and the lowest evaluated merit function is plotted at each iteration. The optimizers are compared by looking at the number of iterations required to reduce the merit function below $10^{-5}$.

% ------------------------------------------
% New section
% ------------------------------------------
\section{Constrained studies for the Bayesian optimizer} \label{Sec_ConOptz_Studies}

% ------------------------------------------
% New subsection
% ------------------------------------------
\subsection{Probability of feasibility and $\ell_2$ penalty} \label{Sec_ConOptz_Studies_ExistingMtds}

\begin{figure}[t!]
	\centering
	\begin{subfigure}[t]{0.328\textwidth}
		\includegraphics[width=\textwidth]{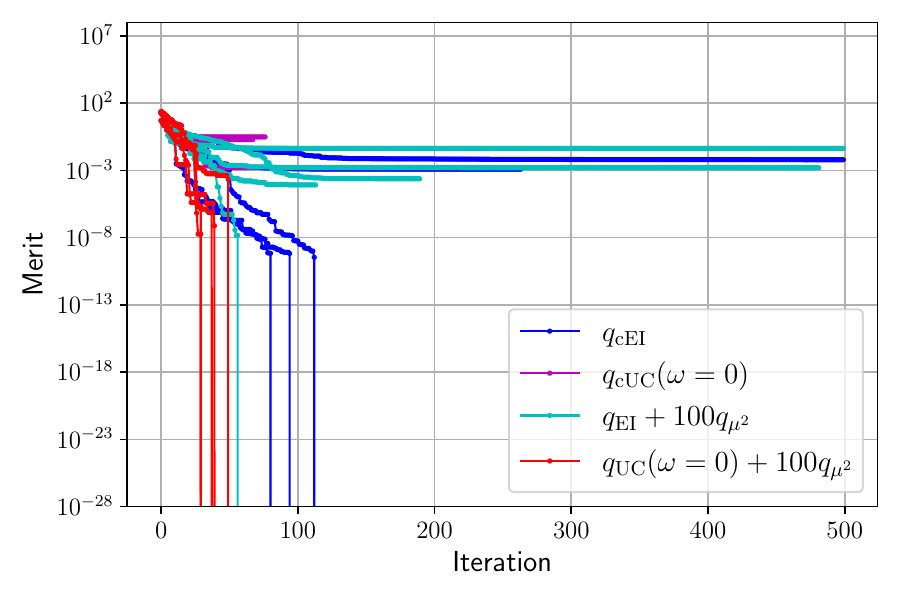}
		\captionsetup{justification=raggedright,singlelinecheck=false}
		\caption{\Eq{Eq_min_quad_nlc_g} with $n_d = 5$: min $f_{\text{quad}}(\xvec)$ with $\xvec^\top \xvec \geq 4$.}
		\label{Fig_ConOptz_MtdExist_fquad_nlc_g_d5}
	\end{subfigure}	
	\begin{subfigure}[t]{0.328\textwidth}
		\includegraphics[width=\textwidth]{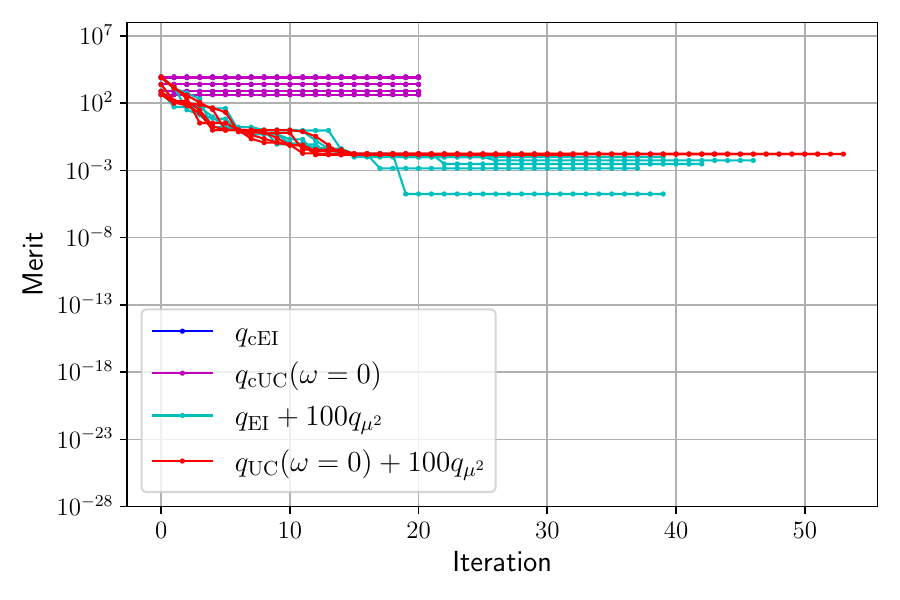}
		\captionsetup{justification=raggedright,singlelinecheck=false}
		\caption{\Eq{Eq_min_prod_nlc_h} with $n_d = 5$: min $f_{\text{prod}}(\xvec)$ with $\xvec^\top \xvec = 1$.}
		\label{Fig_ConOptz_MtdExist_fprod_nlc_h_d5}
	\end{subfigure}	
	\begin{subfigure}[t]{0.328\textwidth}
		\includegraphics[width=\textwidth]{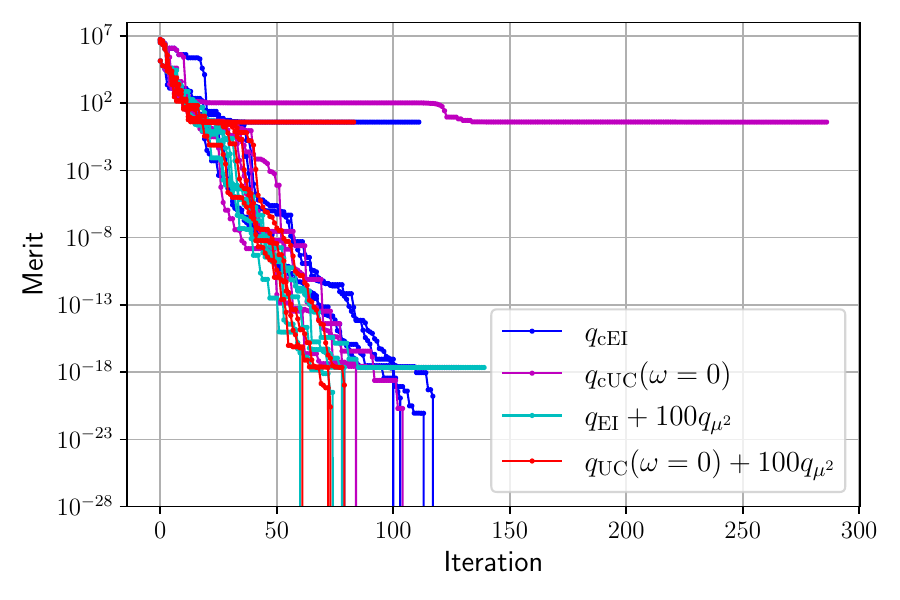}
		\captionsetup{justification=raggedright,singlelinecheck=false}
		\caption{\Eq{Eq_min_nlc_g_Rosen} with $n_d = 5$: min $f_{\text{Rosen}}(\xvec)$ with $\xvec^\top \xvec \leq n_d$.}
		\label{Fig_ConOptz_MtdExist_fRosen_nlc_g_d5}
	\end{subfigure}	
	\begin{subfigure}[t]{0.328\textwidth}
		\includegraphics[width=\textwidth]{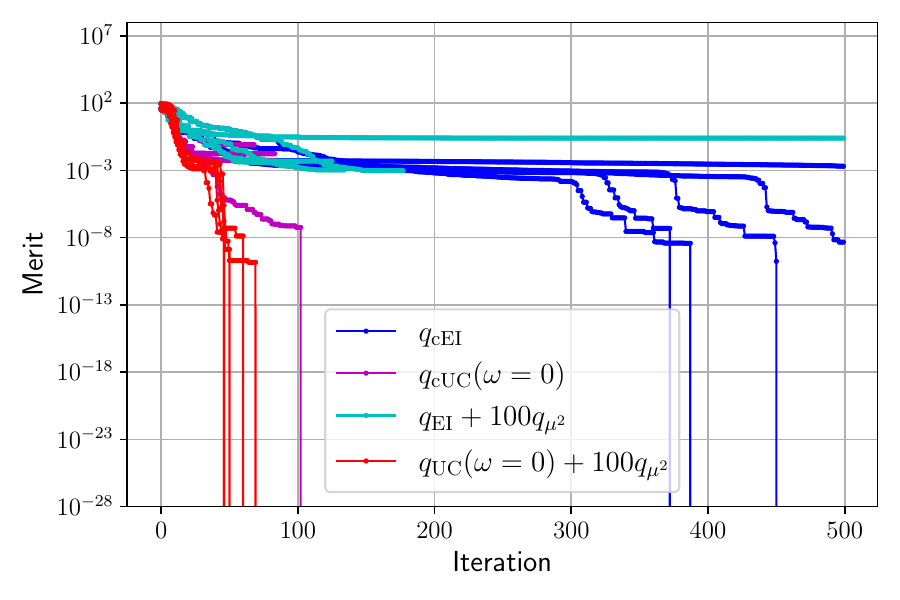}
		\captionsetup{justification=raggedright,singlelinecheck=false}
		\caption{\Eq{Eq_min_quad_nlc_g} with $n_d = 20$: min $f_{\text{quad}}(\xvec)$ with $\xvec^\top \xvec \geq 4$.}
		\label{Fig_ConOptz_MtdExist_fquad_nlc_g_d20}
	\end{subfigure}	
	\begin{subfigure}[t]{0.328\textwidth}
		\includegraphics[width=\textwidth]{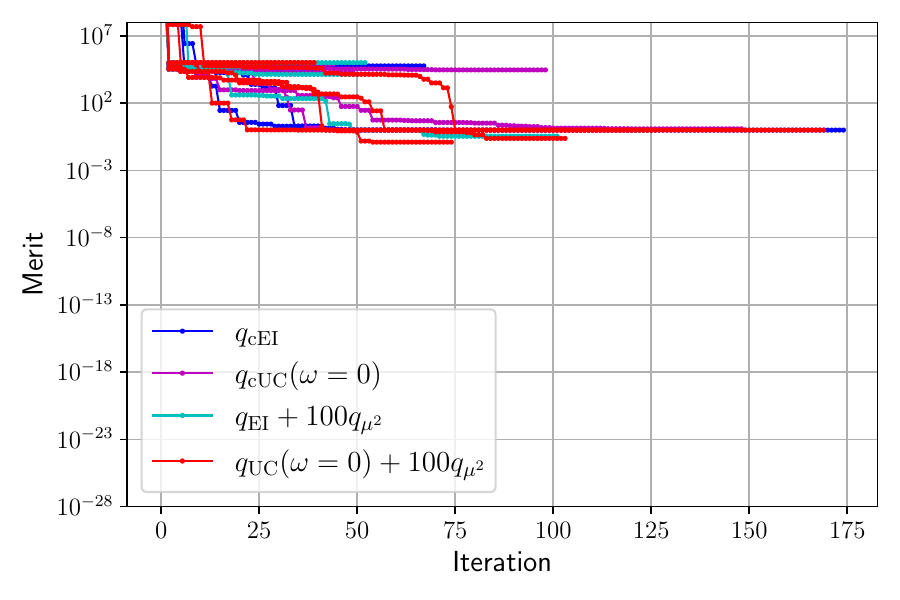}
		\captionsetup{justification=raggedright,singlelinecheck=false}
		\caption{\Eq{Eq_min_prod_nlc_h} with $n_d = 20$: min $f_{\text{prod}}(\xvec)$ with $\xvec^\top \xvec = 1$.}
		\label{Fig_ConOptz_MtdExist_fprod_nlc_h_d20}
	\end{subfigure}	
	\begin{subfigure}[t]{0.328\textwidth}
		\includegraphics[width=\textwidth]{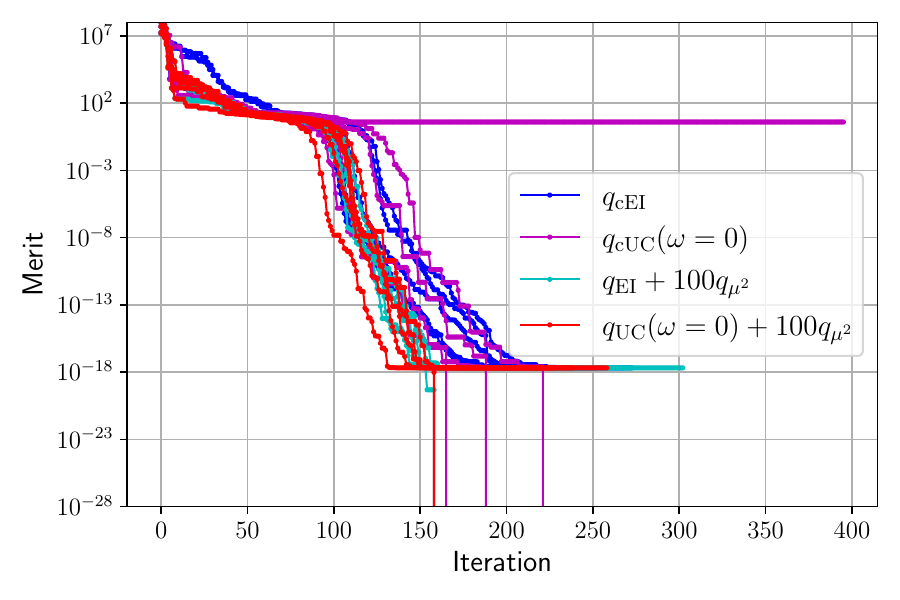}
		\captionsetup{justification=raggedright,singlelinecheck=false}
		\caption{\Eq{Eq_min_nlc_g_Rosen} with $n_d = 20$: min $f_{\text{Rosen}}(\xvec)$ with $\xvec^\top \xvec \leq n_d$.}
		\label{Fig_ConOptz_MtdExist_fRosen_nlc_g_d20}
	\end{subfigure}	
	\caption{Comparison of existing methods for constrained Bayesian optimization presented in \Sec{Sec_ConOptz_MethodsExist} to solve the three constrained test cases from \Sec{Sec_ConOptz_CstrTestCases}. The merit function is from \Eq{Eq_Lag_exa_aug} with $\rho = 100$.}
	\label{Fig_ConOptz_MtdExist}
\end{figure}

Previously developed methods of performing constrained Bayesian optimization are considered in this subsection. In \Sec{Sec_ConOptz_CurrentCstrBo} an acquisition function based on the $\ell_2$ penalty method was introduced along with the probability of feasibility, which was combined with the expected improvement and upper confidence acquisition functions from \Eqs{Eq_acq_cEI}{Eq_acq_cUC}, respectively. These methods of enabling the Bayesian optimizer to solve nonlinearly-constrained optimization problems are applied to the constrained test cases presented in \Sec{Sec_ConOptz_CstrTestCases}. \Fig{Fig_ConOptz_MtdExist} compares the lowest merit function after each evaluation, which comes from the exact augmented Lagrangian from \Eq{Eq_Lag_exa_aug} with $\rho = 100$. For the constrained minimization of $f_{\text{quad}}(\xvec)$ from \Eq{Eq_min_quad_nlc_g}, which is shown in \Figs{Fig_ConOptz_MtdExist_fquad_nlc_g_d5}{Fig_ConOptz_MtdExist_fquad_nlc_g_d20} for $n_d=5$ and $n_d=20$, respectively, the most effective acquisition function combines the $\ell_2$ penalty with the upper confidence acquisition function $\qUC$. The use of the $\ell_2$ penalty with $\qEI$ results in the Bayesian optimizer being unable to reduce the merit function below $10^{-3}$ for all but one of the optimization runs. Similarly, the use of the probability of feasibility for $\qCEI(\xvec)$ and $\qCUC(\xvec; \omega=0)$ results in the optimizer stalling with a high merit function evaluation or taking significantly more iterations to reduce the merit function below $10^{-5}$ than with the use of the upper confidence and $\ell_2$ penalty acquisition functions. The Bayesian optimizer with any of the four different acquisition functions struggles to handle the nonlinear equality constraint in \Figs{Fig_ConOptz_MtdExist_fprod_nlc_h_d5}{Fig_ConOptz_MtdExist_fprod_nlc_h_d20} for $n_d=5$ and $n_d=20$, respectively. 

%--------------------------------------------------------------------------
%-------------------------------------------------------------------------- 

% ------------------------------------------
% New subsection
% ------------------------------------------
\subsection{Exact augmented Lagrangian with Bayesian optimization} \label{Sec_ConOptz_Studies_ExaLag}

\begin{figure}[t!]
	\centering
	\begin{subfigure}[t]{0.328\textwidth}
		\includegraphics[width=\textwidth]{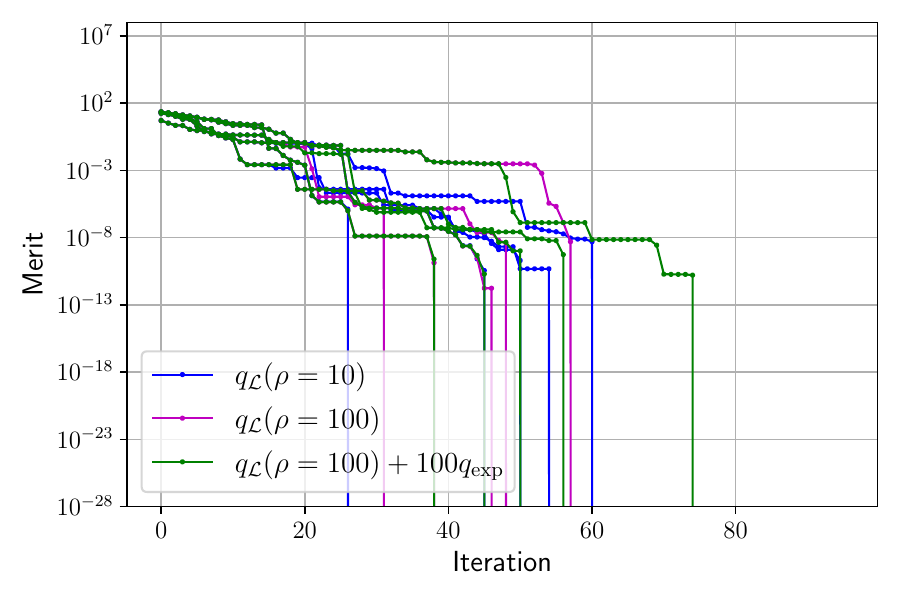}
		\captionsetup{justification=raggedright,singlelinecheck=false}
		\caption{\Eq{Eq_min_quad_nlc_g} with $n_d = 5$: min $f_{\text{quad}}(\xvec)$ with $\xvec^\top \xvec \geq 4$.}
		\label{Fig_ConOptz_ExaLag_fquad_nlc_g_d5}
	\end{subfigure}	
	\begin{subfigure}[t]{0.328\textwidth}
		\includegraphics[width=\textwidth]{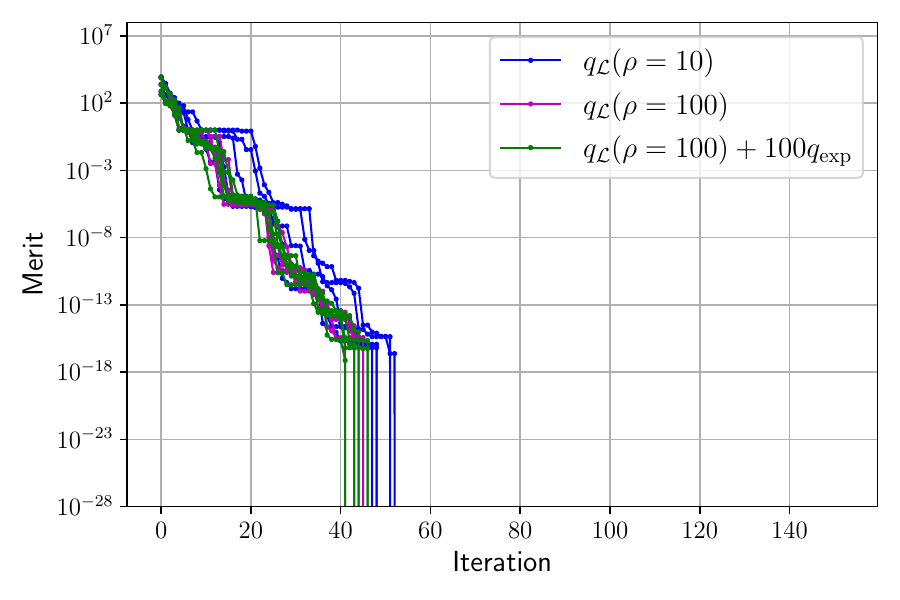}
		\captionsetup{justification=raggedright,singlelinecheck=false}
		\caption{\Eq{Eq_min_prod_nlc_h} with $n_d = 5$: min $f_{\text{prod}}(\xvec)$ with $\xvec^\top \xvec = 1$.}
		\label{Fig_ConOptz_ExaLag_fprod_nlc_h_d5}
	\end{subfigure}	
	\begin{subfigure}[t]{0.328\textwidth}
		\includegraphics[width=\textwidth]{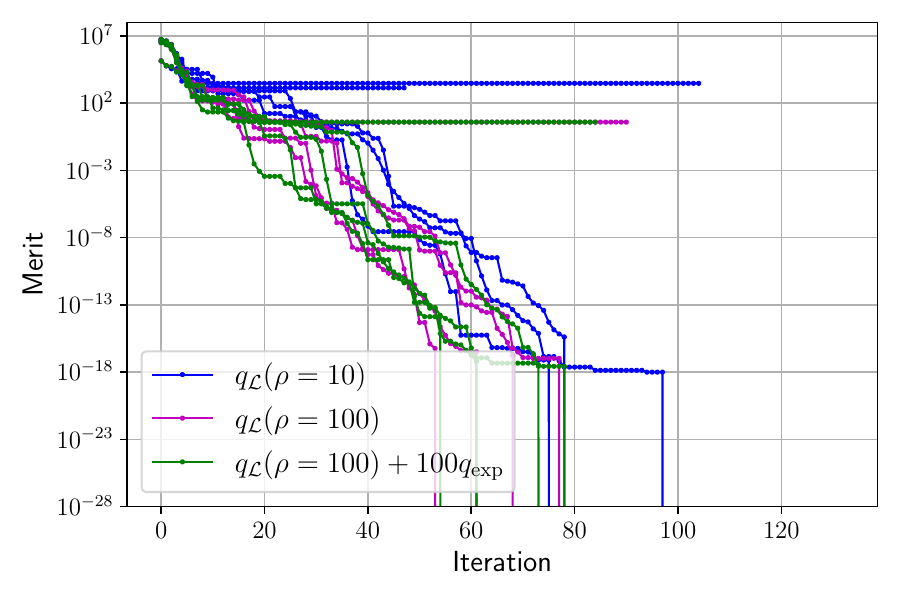}
		\captionsetup{justification=raggedright,singlelinecheck=false}
		\caption{\Eq{Eq_min_nlc_g_Rosen} with $n_d = 5$: min $f_{\text{Rosen}}(\xvec)$ with $\xvec^\top \xvec \leq n_d$.}
		\label{Fig_ConOptz_ExaLag_fRosen_nlc_g_d5}
	\end{subfigure}	
	\begin{subfigure}[t]{0.328\textwidth}
		\includegraphics[width=\textwidth]{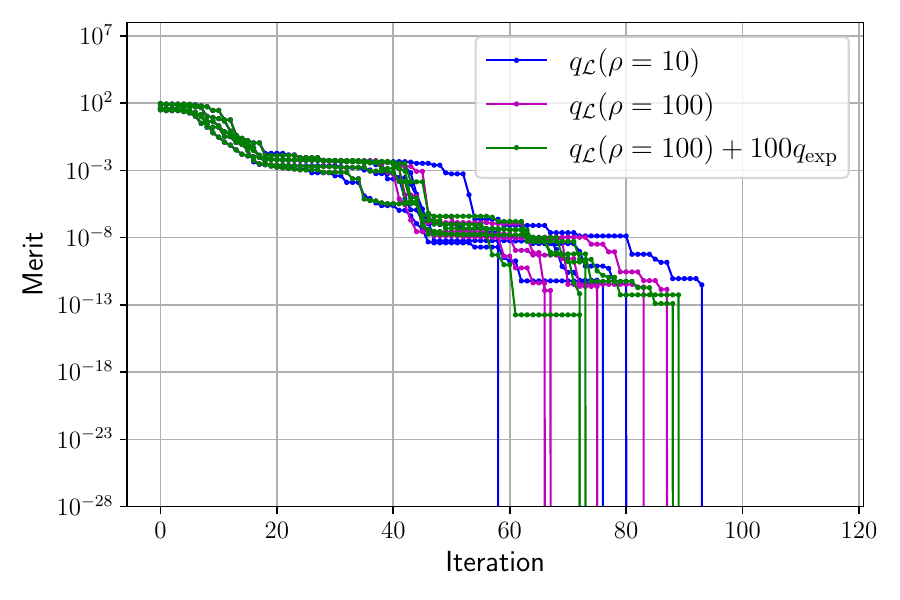}
		\captionsetup{justification=raggedright,singlelinecheck=false}
		\caption{\Eq{Eq_min_quad_nlc_g} with $n_d = 20$: min $f_{\text{quad}}(\xvec)$ with $\xvec^\top \xvec \geq 4$.}
		\label{Fig_ConOptz_ExaLag_fquad_nlc_g_d20}
	\end{subfigure}	
	\begin{subfigure}[t]{0.328\textwidth}
		\includegraphics[width=\textwidth]{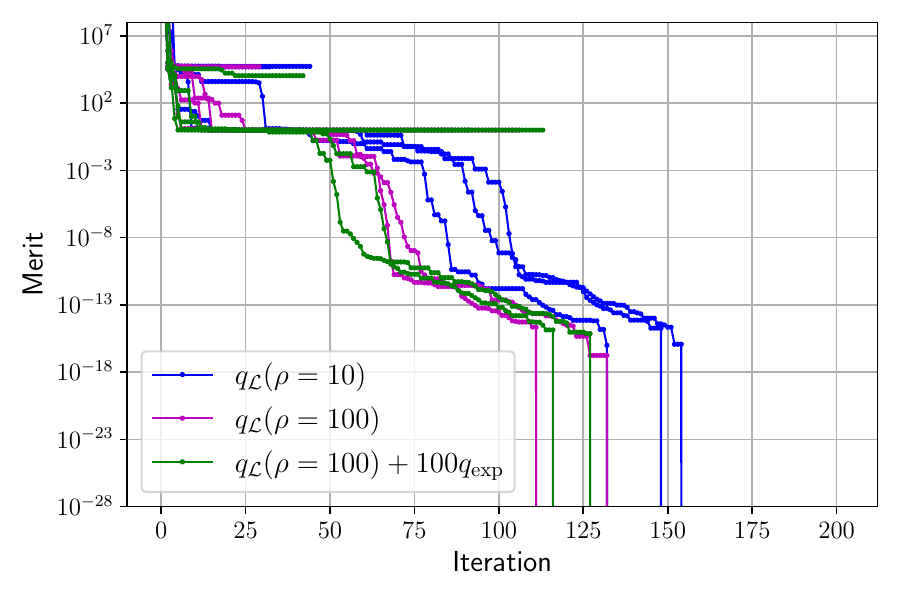}
		\captionsetup{justification=raggedright,singlelinecheck=false}
		\caption{\Eq{Eq_min_prod_nlc_h} with $n_d = 20$: min $f_{\text{prod}}(\xvec)$ with $\xvec^\top \xvec = 1$.}
		\label{Fig_ConOptz_ExaLag_fprod_nlc_h_d20}
	\end{subfigure}	
	\begin{subfigure}[t]{0.328\textwidth}
		\includegraphics[width=\textwidth]{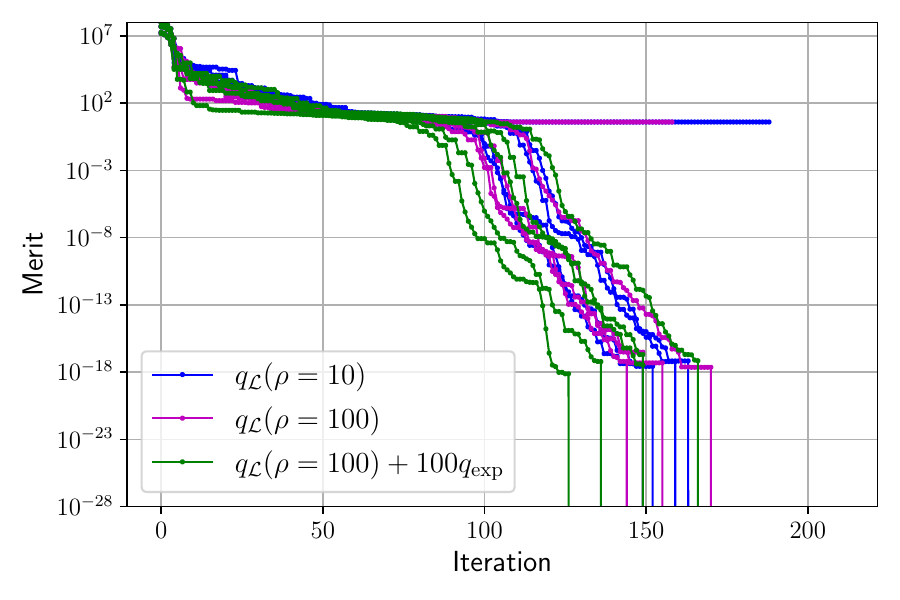}
		\captionsetup{justification=raggedright,singlelinecheck=false}
		\caption{\Eq{Eq_min_nlc_g_Rosen} with $n_d = 20$: min $f_{\text{Rosen}}(\xvec)$ with $\xvec^\top \xvec \leq n_d$.}
		\label{Fig_ConOptz_ExaLag_fRosen_nlc_g_d20}
	\end{subfigure}	
	\caption[Study with the exact augmented Lagrangian for constrained Bayesian optimization.]{The exact augmented Lagrangian acquisition function from \Sec{Sec_ConOptz_MethodsNew_ExactAugLag} is compared with different settings to solve the constrained test cases from \Sec{Sec_ConOptz_CstrTestCases}. \Eq{Eq_Lag_exa_aug} with $\rho = 100$ is used as the merit function.}
	\label{Fig_ConOptz_ExaLag}
\end{figure}

The exact augmented Lagrangian acquisition function $\qLag(\xvec; \rho)$ from \Eq{Eq_acq_Lag_exa_aug} is applied in this subsection with different values of $\rho$, which is the coefficient for an $\ell_2$ penalty, along with the probabilistic acquisition function $\qExp(\xvec)$ from \Eq{Eq_acq_qExp}. The merit function from \Eq{Eq_Lag_exa_aug} with $\rho = 100$ is plotted in \Fig{Fig_ConOptz_ExaLag} for the three constrained test cases from \Sec{Sec_ConOptz_MethodsNew_ProbPenalty}. For all three constrained test cases with both $n_d=5$ and $n_d=20$, the three variations of the exact augmented Lagrangian acquisition function enable the Bayesian optimizer to converge the merit function below $10^{-5}$ for nearly all optimization runs. This demonstrates that the exact augmented Lagrangian acquisition function is effective for both nonlinear inequality and equality constraints.

%--------------------------------------------------------------------------
%-------------------------------------------------------------------------- 

\begin{figure}[t!]
	\centering
	\begin{subfigure}[t]{0.328\textwidth}
		\includegraphics[width=\textwidth]{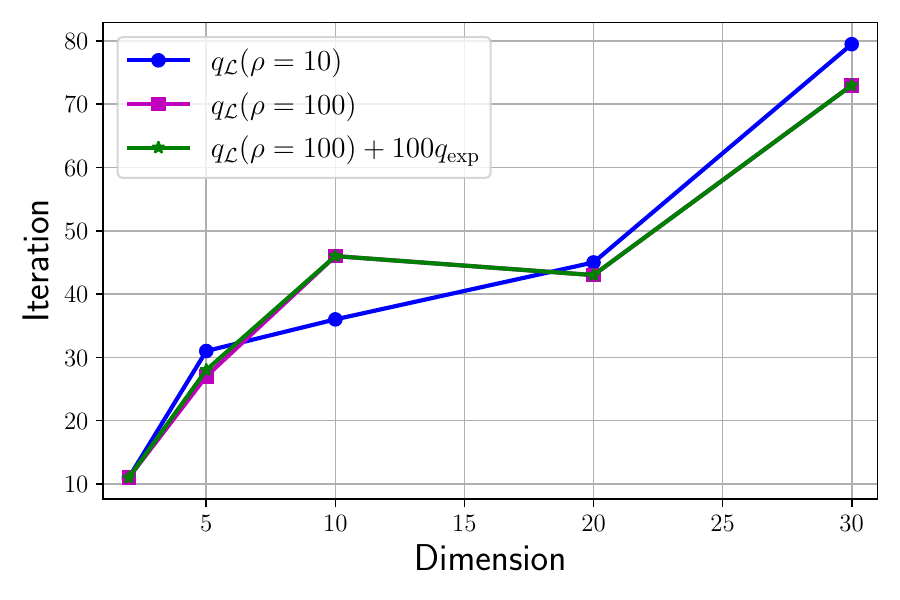}
		\captionsetup{justification=raggedright,singlelinecheck=false}
		\caption{min $f_{\text{quad}}(\xvec)$ with $\xvec^\top \xvec \geq 4$ from \Eq{Eq_min_quad_nlc_g}.}
		\label{Fig_ConOptz_ExaLag_tol_fquad}
	\end{subfigure}	
	\begin{subfigure}[t]{0.328\textwidth}
		\includegraphics[width=\textwidth]{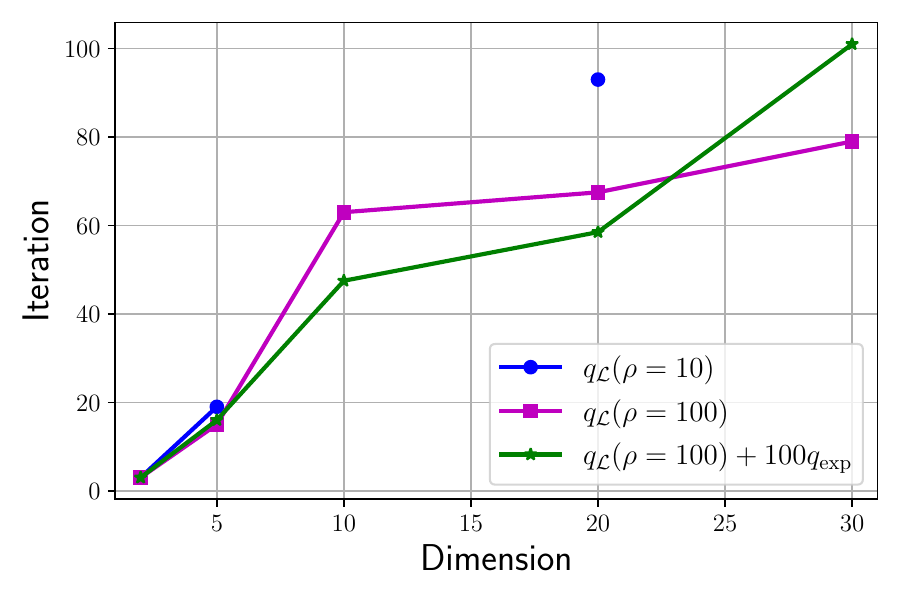}
		\captionsetup{justification=raggedright,singlelinecheck=false}
		\caption{min $f_{\text{prod}}(\xvec)$ with $\xvec^\top \xvec = 1$ from \Eq{Eq_min_prod_nlc_h}.}
		\label{Fig_ConOptz_ExaLag_tol_fprod}
	\end{subfigure}	
	\begin{subfigure}[t]{0.328\textwidth}
		\includegraphics[width=\textwidth]{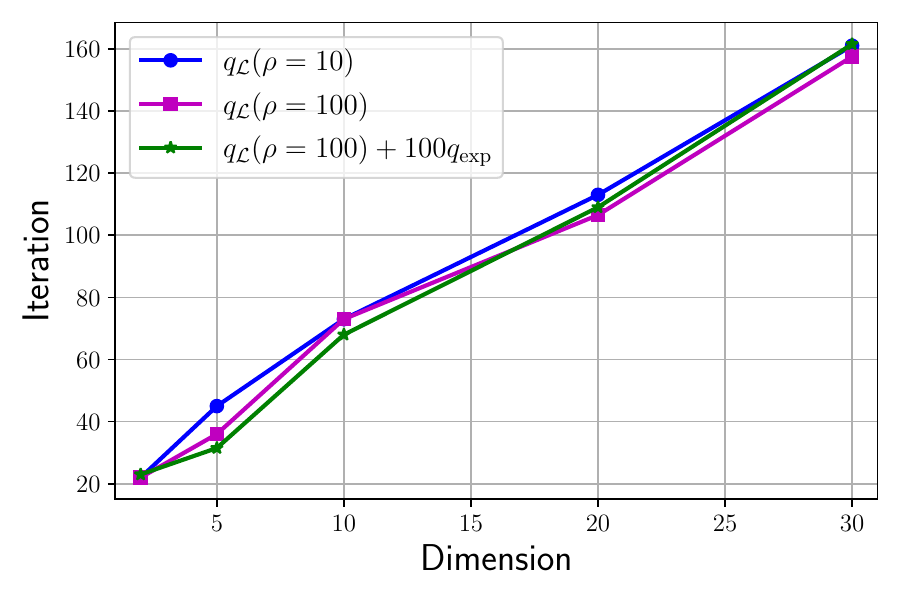}
		\captionsetup{justification=raggedright,singlelinecheck=false}
		\caption{min $f_{\text{Rosen}}(\xvec)$ with $\xvec^\top \xvec \leq n_d$ from \Eq{Eq_min_nlc_g_Rosen}.}
		\label{Fig_ConOptz_ExaLag_tol_fRosen}
	\end{subfigure}	
	\caption[Median number of iterations for the Bayesian optimizer using an exact augmented Lagrangian to reduce the merit function below $10^{-5}$.]{Median number of iterations for the Bayesian optimizer using the exact augmented Lagrangian $\qLag(\xvec)$ from \Eq{Eq_acq_Lag_exa_aug} with different penalty functions to reduce the merit function from \Eq{Eq_Lag_exa_aug} with $\rho = 100$ below $10^{-5}$.}
	\label{Fig_ConOptz_ExaLag_tol}
\end{figure}

\Fig{Fig_ConOptz_ExaLag_tol} shows the number of iterations required for the Bayesian optimizer using the exact augmented Lagrangian with different penalty terms to reduce the merit function below $10^{-5}$. The largest difference in the results lies in \Fig{Fig_ConOptz_ExaLag_tol_fprod}, which is for the test case with an equality constraint. It is evident that increasing $\rho$, which is the coefficient for the $\ell_2$ penalty, from 10 to 100 is beneficial in this case. The addition of $\qExp(\xvec)$ to $\qLag(\xvec)$ only has a noticeable impact on the test case with the equality constraint shown in \Fig{Fig_ConOptz_ExaLag_tol_fprod}. In this case, the addition of this acquisition function leads to a reduction in the number of iterations required to reach the desired tolerance when $n_d = 10, 20$, an increase when $n_d=30$, and minor differences for $n_d =2,5$. For the three test cases considered, the constraints all involved $\xvec^\top \xvec$, which is a convex function. Therefore, the gradient for each of these constraints always points in the direction where they are satisfied. However, for more complicated nonlinear constraints, this may not always be the case, which is when the additional exploration provided by $\qExp(\xvec)$ could prove helpful. Additional test cases are needed to better understand the impact of $\qExp(\xvec)$ on the performance of the Bayesian optimizer.

%--------------------------------------------------------------------------
%-------------------------------------------------------------------------- 

% ------------------------------------------
% New subsection
% ------------------------------------------
\subsection{Strong enforcement of the constraints} \label{Sec_ConOptz_Studies_Strong}

\begin{figure}[t!]
	\centering
	\begin{subfigure}[t]{0.328\textwidth}
		\includegraphics[width=\textwidth]{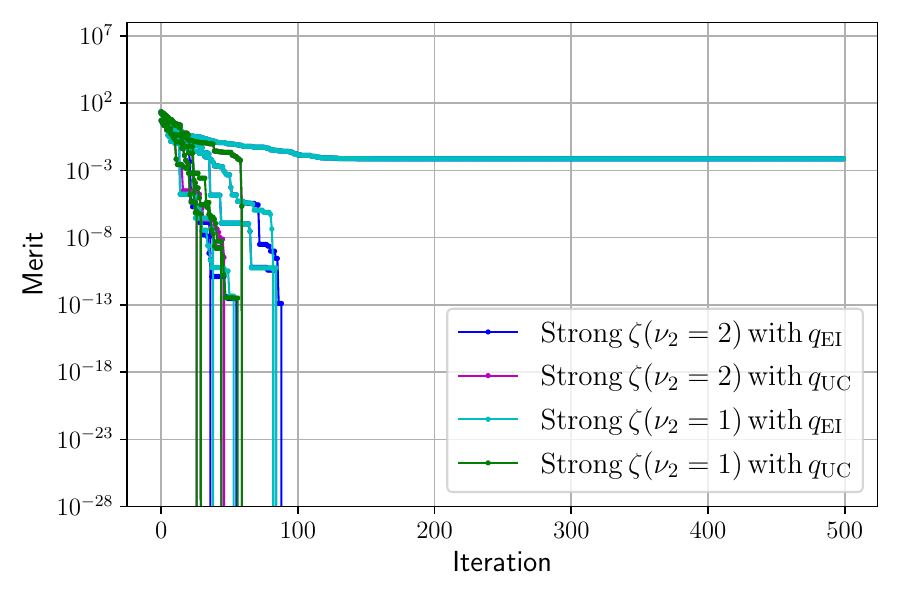}
		\captionsetup{justification=raggedright,singlelinecheck=false}
		\caption{\Eq{Eq_min_quad_nlc_g} with $n_d = 5$: min $f_{\text{quad}}(\xvec)$ with $\xvec^\top \xvec \geq 4$.}
		\label{Fig_ConOptz_StudyStrong_fquad_nlc_g_d5}
	\end{subfigure}	
	\begin{subfigure}[t]{0.328\textwidth}
		\includegraphics[width=\textwidth]{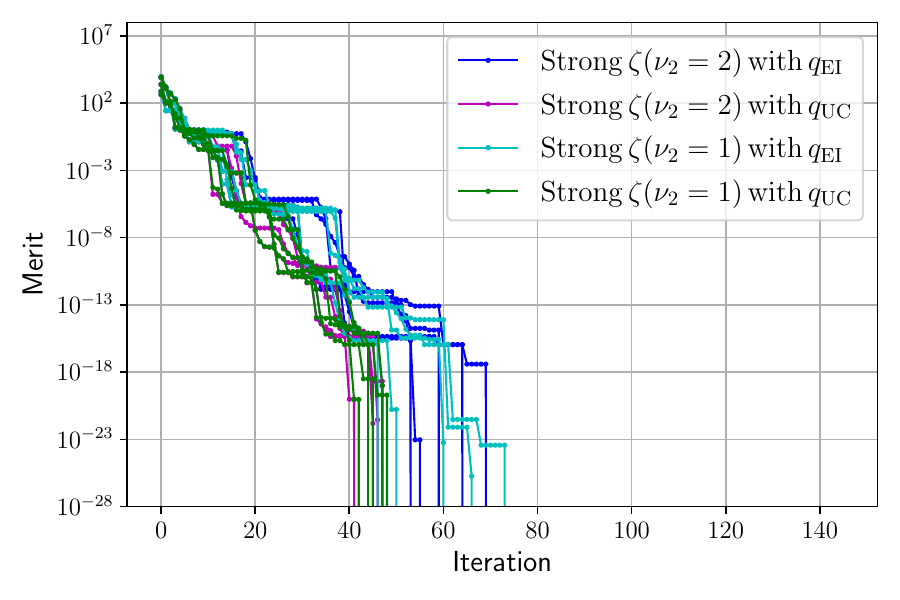}
		\captionsetup{justification=raggedright,singlelinecheck=false}
		\caption{\Eq{Eq_min_prod_nlc_h} with $n_d = 5$: min $f_{\text{prod}}(\xvec)$ with $\xvec^\top \xvec = 1$.}
		\label{Fig_ConOptz_StudyStrong_fprod_nlc_h_d5}
	\end{subfigure}	
	\begin{subfigure}[t]{0.328\textwidth}
		\includegraphics[width=\textwidth]{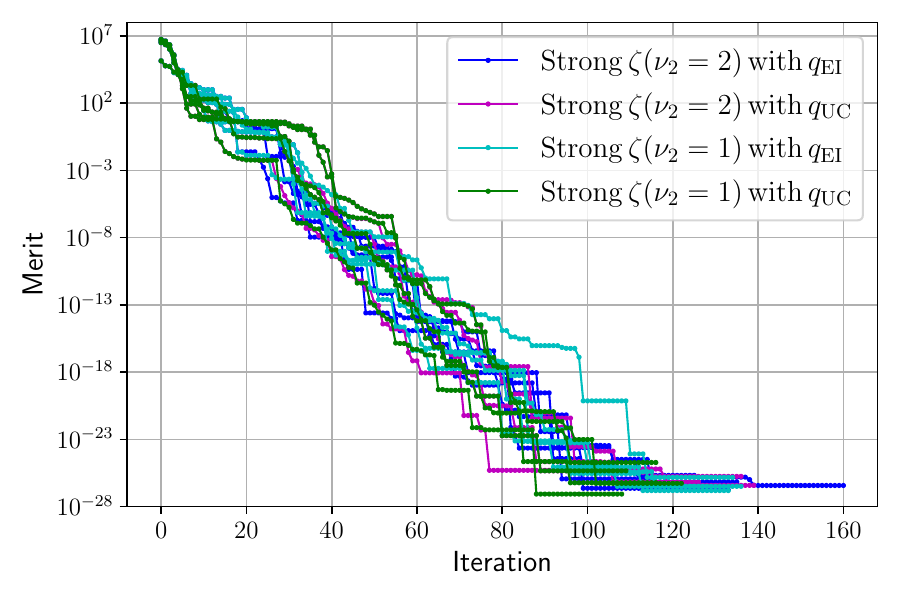}
		\captionsetup{justification=raggedright,singlelinecheck=false}
		\caption{\Eq{Eq_min_nlc_g_Rosen} with $n_d = 5$: min $f_{\text{Rosen}}(\xvec)$ with $\xvec^\top \xvec \leq n_d$.}
		\label{Fig_ConOptz_StudyStrong_fRosen_nlc_g_d5}
	\end{subfigure}	
	\begin{subfigure}[t]{0.328\textwidth}
		\includegraphics[width=\textwidth]{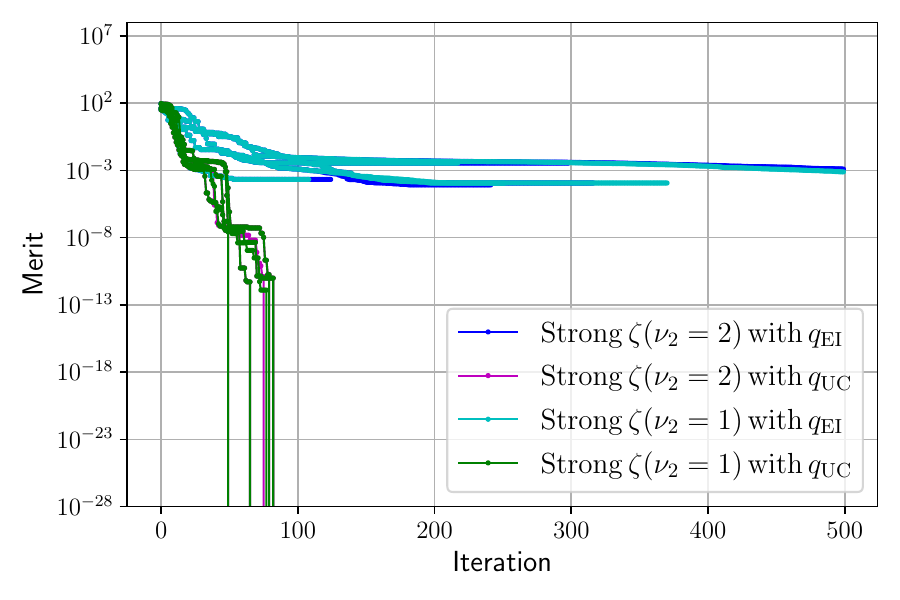}
		\captionsetup{justification=raggedright,singlelinecheck=false}
		\caption{\Eq{Eq_min_quad_nlc_g} with $n_d = 20$: min $f_{\text{quad}}(\xvec)$ with $\xvec^\top \xvec \geq 4$.}
		\label{Fig_ConOptz_StudyStrong_fquad_nlc_g_d20}
	\end{subfigure}	
	\begin{subfigure}[t]{0.328\textwidth}
		\includegraphics[width=\textwidth]{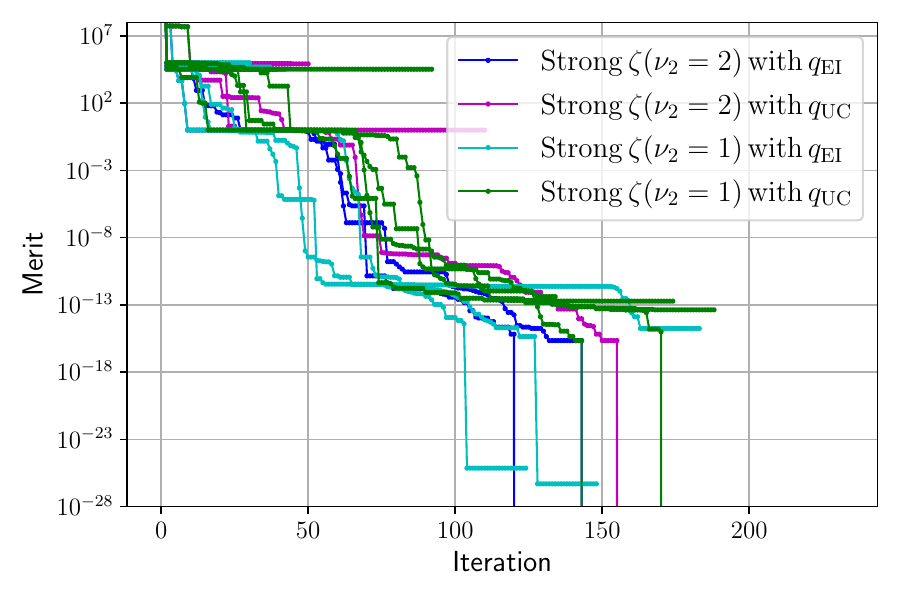}
		\captionsetup{justification=raggedright,singlelinecheck=false}
		\caption{\Eq{Eq_min_prod_nlc_h} with $n_d = 20$: min $f_{\text{prod}}(\xvec)$ with $\xvec^\top \xvec = 1$.}
		\label{Fig_ConOptz_StudyStrong_fprod_nlc_h_d20}
	\end{subfigure}	
	\begin{subfigure}[t]{0.328\textwidth}
		\includegraphics[width=\textwidth]{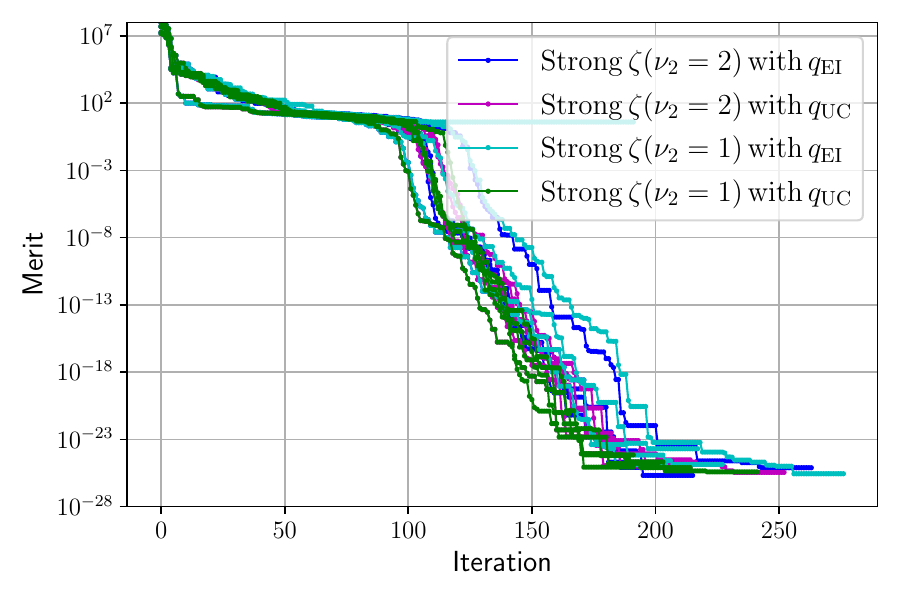}
		\captionsetup{justification=raggedright,singlelinecheck=false}
		\caption{\Eq{Eq_min_nlc_g_Rosen} with $n_d = 20$: min $f_{\text{Rosen}}(\xvec)$ with $\xvec^\top \xvec \leq n_d$.}
		\label{Fig_ConOptz_StudyStrong_fRosen_nlc_g_d20}
	\end{subfigure}	
	\caption[Study with strong enforcement of the nonlinear constraint for the Bayesian optimizer.]{The test cases from \Sec{Sec_ConOptz_CstrTestCases} are solved using the Bayesian optimizer with strong enforcement of the nonlinear constraints, as detailed in Algorithm~\ref{Alg_StrongEnf}. The results are compared with the merit function coming from \Eq{Eq_Lag_exa_aug} with $\rho = 100$ and $\omega = 0$ is used for $\qUC$ from \Eq{Eq_acq_UC}.}
	\label{Fig_ConOptz_StudyStrong}
\end{figure}

The strong enforcement of the nonlinear constraints from \Sec{Sec_ConOptz_MethodsNew_StrongEnf} is investigated in this section using the constrained test cases from \Sec{Sec_ConOptz_CstrTestCases}. The upper confidence and expected improvement acquisition function from \Eqs{Eq_acq_UC}{Eq_acq_EI}, respectively, are used and Algorithm~\ref{Alg_StrongEnf} is used to enforce the nonlinear constraints with the sigmoid function $\zeta(z; \nu_1 = 10, \nu_2)$ from \Eq{Eq_sigmoid_strong_enf} with either $\nu_2 = 1$ or $\nu_2 = 2$. It is clear from \Fig{Fig_ConOptz_StudyStrong} that the results for the Bayesian optimizer using the four combinations of acquisition functions are similar in nearly all cases. The exception is for the constrained minimization of $f_{\text{quad}}(\xvec)$ from \Eq{Eq_min_quad_nlc_g} for $n_d=20$, which is shown in \Fig{Fig_ConOptz_StudyStrong_fquad_nlc_g_d20}. In this case, the Bayesian optimizer that uses the expected improvement acquisition function with either $\nu_2 = 1$ or $\nu_2=2$ does not converge the merit function below $10^{-4}$ for any of the five independent optimization runs. The expected improvement acquisition function promotes more exploration than the upper confidence acquisition function, which is clearly impacting the convergence of the Bayesian optimizer for the constrained minimization of $f_{\text{quad}}(\xvec)$ from \Eq{Eq_min_quad_nlc_g}. This same trend was also observed in \Sec{Sec_ConOptz_Studies_ExistingMtds} when the expected improvement acquisition function was used with the $\ell_2$ penalty and the probability of feasibility.

%--------------------------------------------------------------------------
%-------------------------------------------------------------------------- 

\begin{figure}[t!]
	\centering
	\begin{subfigure}[t]{0.328\textwidth}
		\includegraphics[width=\textwidth]{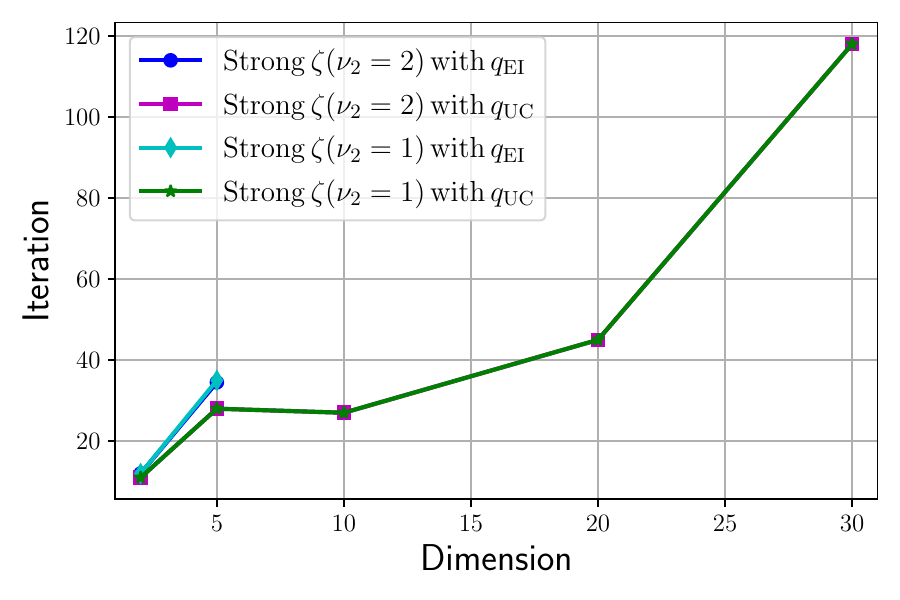}
		\captionsetup{justification=raggedright,singlelinecheck=false}
		\caption{min $f_{\text{quad}}(\xvec)$ with $\xvec^\top \xvec \geq 4$ from \Eq{Eq_min_quad_nlc_g}.}
		\label{Fig_ConOptz_StrongEnf_tol_fquad}
	\end{subfigure}	
	\begin{subfigure}[t]{0.328\textwidth}
		\includegraphics[width=\textwidth]{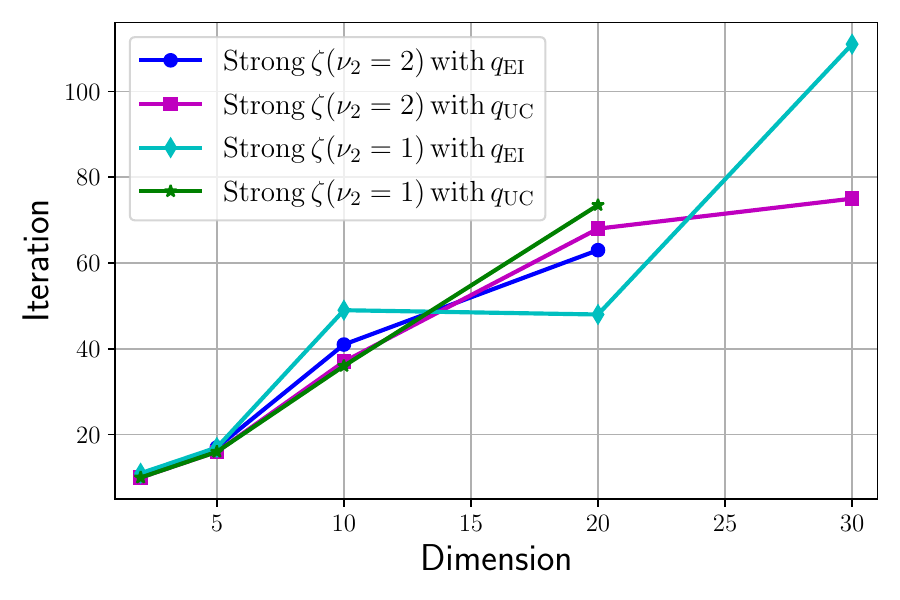}
		\captionsetup{justification=raggedright,singlelinecheck=false}
		\caption{min $f_{\text{prod}}(\xvec)$ with $\xvec^\top \xvec = 1$ from \Eq{Eq_min_prod_nlc_h}.}
		\label{Fig_ConOptz_StrongEnf_tol_fprod}
	\end{subfigure}	
	\begin{subfigure}[t]{0.328\textwidth}
		\includegraphics[width=\textwidth]{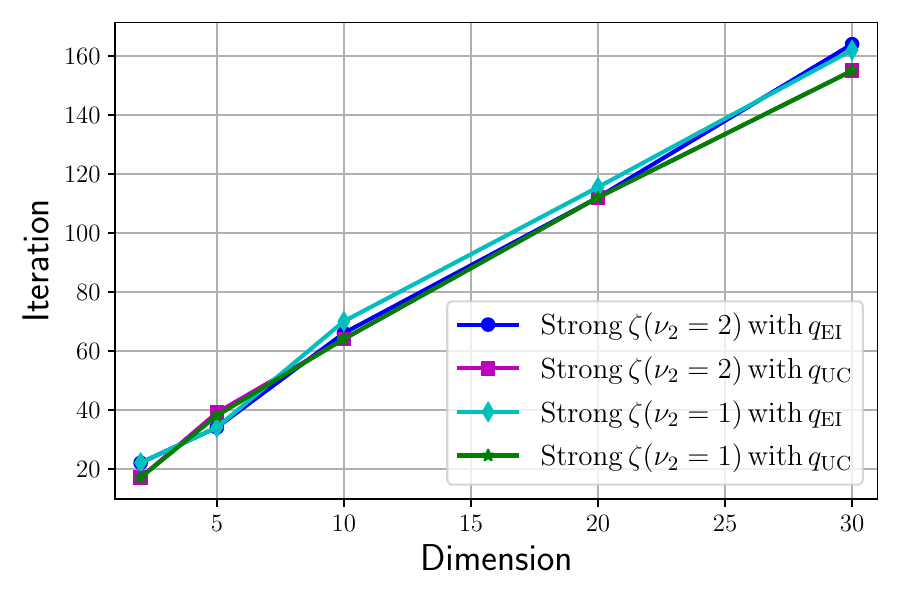}
		\captionsetup{justification=raggedright,singlelinecheck=false}
		\caption{min $f_{\text{Rosen}}(\xvec)$ with $\xvec^\top \xvec \leq n_d$ from \Eq{Eq_min_nlc_g_Rosen}.}
		\label{Fig_ConOptz_StrongEnf_tol_fRosen}
	\end{subfigure}	
	\caption[Median number of iterations for the Bayesian optimizer using a strong enforcement of the nonlinear constraints to reduce the merit function below $10^{-5}$.]{The median number of iterations for the Bayesian optimizer using Algorithm~\ref{Alg_StrongEnf} for strong enforcement of the nonlinear constraints to reduce the merit function below $10^{-5}$. For the acquisition $\qUC$ from \Eq{Eq_acq_UC} the hyerparameter $\omega = 0$ was used and the merit function is the exact augmented Lagrangian from \Eq{Eq_Lag_exa_aug} with $\rho = 100$.}
	\label{Fig_ConOptz_StrongEnf_tol}
\end{figure}

In \Fig{Fig_ConOptz_StrongEnf_tol} the median number of iterations required for the Bayesian optimizer to reduce the merit function below $10^{-5}$ is shown for the four combinations of acquisition functions. \Fig{Fig_ConOptz_StrongEnf_tol_fquad} shows that the use of the expected improvement acquisition function does not enable the Bayesian optimizer to reach the desired tolerance for any of the optimization runs for the constrained minimization of $f_{\text{quad}}(\xvec)$ from \Eq{Eq_min_quad_nlc_g} for $n_d \geq 10$. The results for the use of the upper confidence acquisition function with either $\nu_2 = 1$ or $\nu_2 = 2$ are similar for all three test cases. They are nearly identical in \Figs{Fig_ConOptz_StrongEnf_tol_fquad}{Fig_ConOptz_StrongEnf_tol_fRosen} since these test cases have inequality constraints. When the constraint at $\xbest$ is satisfied, the strong enforcement method results in identical minimizations of the acquisition function, regardless of the value of $\nu_2$. The upper confidence acquisition function with $\nu_2=1$ for the sigmoid function from \Eq{Eq_sigmoid_strong_enf} was selected as the default setting when using Algorithm~\ref{Alg_StrongEnf} to solve nonlinearly-constrained optimization problems. 

In the following section, the Bayesian optimizer either with the exact augmented Lagrangian acquisition function or with strong enforcement of the nonlinear constraints is compared to quasi-Newton optimizers from MATLAB and SciPy.

%--------------------------------------------------------------------------
%-------------------------------------------------------------------------- 

%% file: Sections/Sec_BoVsQN.tex
%!TEX root = ../Marchildon_ZinggCstrBayeOptz.tex

% ------------------------------------------
% New section
% ------------------------------------------
\section{Comparing Bayesian and quasi-Newton constrained optimizers} \label{Sec_BoVsQN}

% SciPy trust-constr is a qN optimizer: https://docs.scipy.org/doc/scipy/reference/generated/scipy.optimize.minimize.html

The Bayesian optimizer is compared with MATLAB fmincon SQP, SciPy SLSQP, and SciPy trust-constr, which are all quasi-Newton optimizers. The only default parameters that were changed for these optimizers are their tolerances, which were made smaller to ensure the optimizers converged as deeply as possible. The tolerances ``xtol'' and ``gtol'' were reduced to $10^{-16}$ from $10^{-8}$ for the SciPy trust-constr optimizer. Likewise, for the MATLAB fminunc optimizer the options ``StepTolerance'' and ``OptimalityTolerance'' were reduced to $10^{-16}$ from their default value of $10^{-6}$.

Two variations of the Bayesian optimizer are considered to solve the nonlinearly-constrained test cases from \Sec{Sec_ConOptz_CstrTestCases}. The first uses the acquisition function $\qLag(\xvec)$, which is based on an exact augmented Lagrangian, as detailed in Algorithm~\ref{Alg_ExaAugLag}. The second variation for the Bayesian optimizer uses Algorithm~\ref{Alg_StrongEnf} to strongly enforce the nonlinear constraints. The merit function that is used to compare these optimizers is the exact augmented Lagrangian from \Eq{Eq_Lag_exa_aug} with $\rho = 100$.

\begin{figure}[t!]
	\centering
	\begin{subfigure}[t]{0.328\textwidth}
		\includegraphics[width=\textwidth]{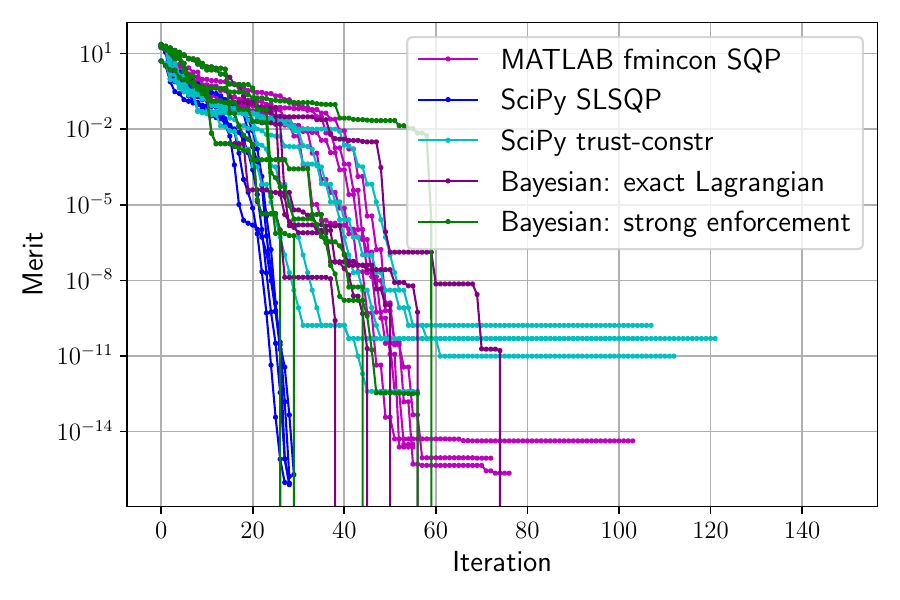}
		\captionsetup{justification=raggedright,singlelinecheck=false}
		\caption{\Eq{Eq_min_quad_nlc_g} with $n_d = 5$: min $f_{\text{quad}}(\xvec)$ with $\xvec^\top \xvec \geq 4$.}
		\label{Fig_ConOptz_qN_vs_BO_fquad_nlc_g_d5}
	\end{subfigure}	
	\begin{subfigure}[t]{0.328\textwidth}
		\includegraphics[width=\textwidth]{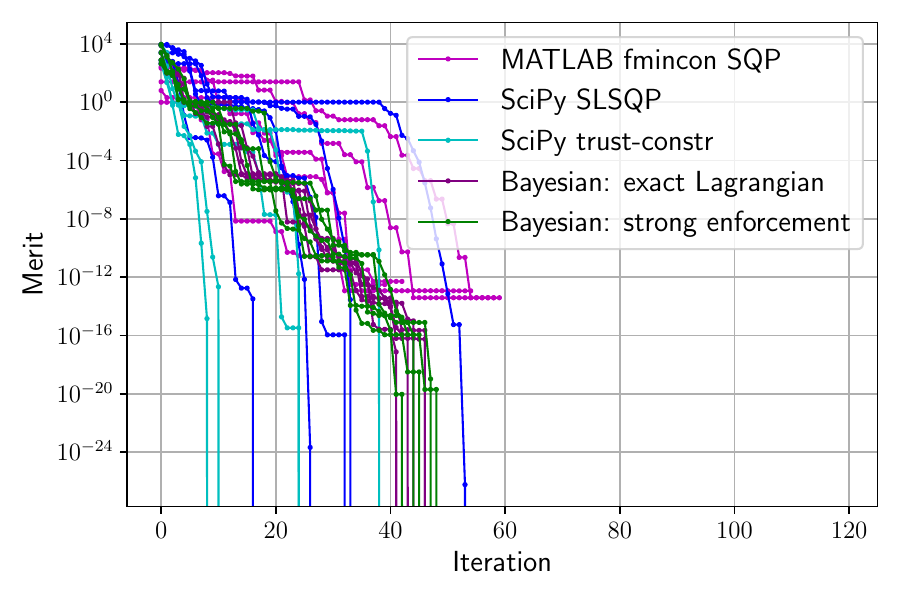}
		\captionsetup{justification=raggedright,singlelinecheck=false}
		\caption{\Eq{Eq_min_prod_nlc_h} with $n_d = 5$: min $f_{\text{prod}}(\xvec)$ with $\xvec^\top \xvec = 1$.}
		\label{Fig_ConOptz_qN_vs_BO_fprod_nlc_h_d5}
	\end{subfigure}	
	\begin{subfigure}[t]{0.328\textwidth}
		\includegraphics[width=\textwidth]{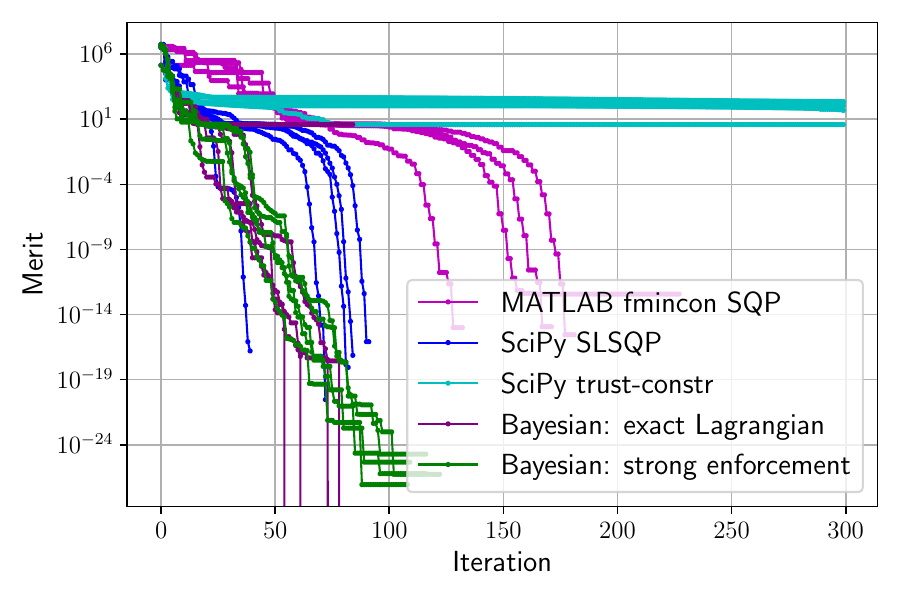}
		\captionsetup{justification=raggedright,singlelinecheck=false}
		\caption{\Eq{Eq_min_nlc_g_Rosen} with $n_d = 5$: min $f_{\text{Rosen}}(\xvec)$ with $\xvec^\top \xvec \leq n_d$.}
		\label{Fig_ConOptz_qN_vs_BO_fRosen_nlc_g_d5}
	\end{subfigure}	
	\begin{subfigure}[t]{0.328\textwidth}
		\includegraphics[width=\textwidth]{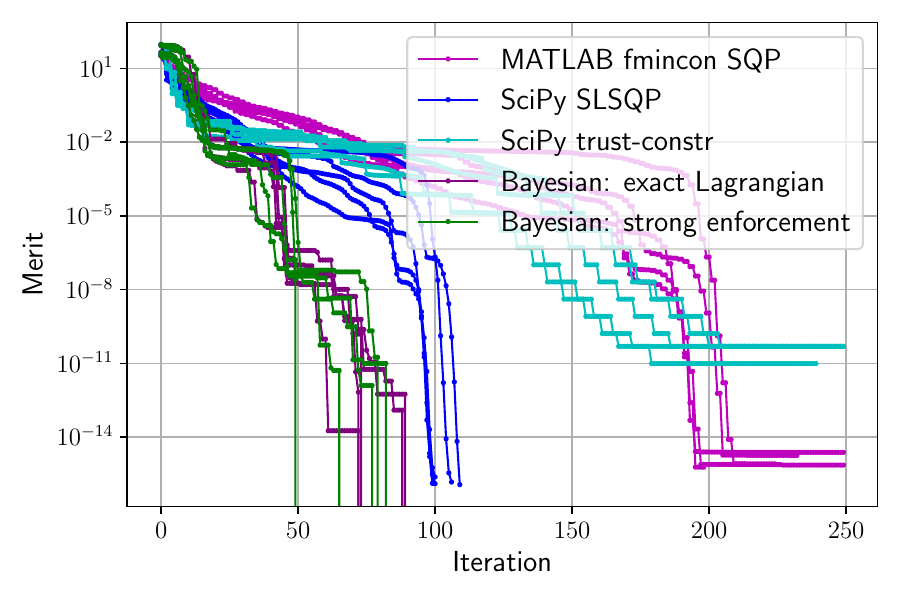}
		\captionsetup{justification=raggedright,singlelinecheck=false}
		\caption{\Eq{Eq_min_quad_nlc_g} with $n_d = 20$: min $f_{\text{quad}}(\xvec)$ with $\xvec^\top \xvec \geq 4$.}
		\label{Fig_ConOptz_qN_vs_BO_fquad_nlc_g_d20}
	\end{subfigure}	
	\begin{subfigure}[t]{0.328\textwidth}
		\includegraphics[width=\textwidth]{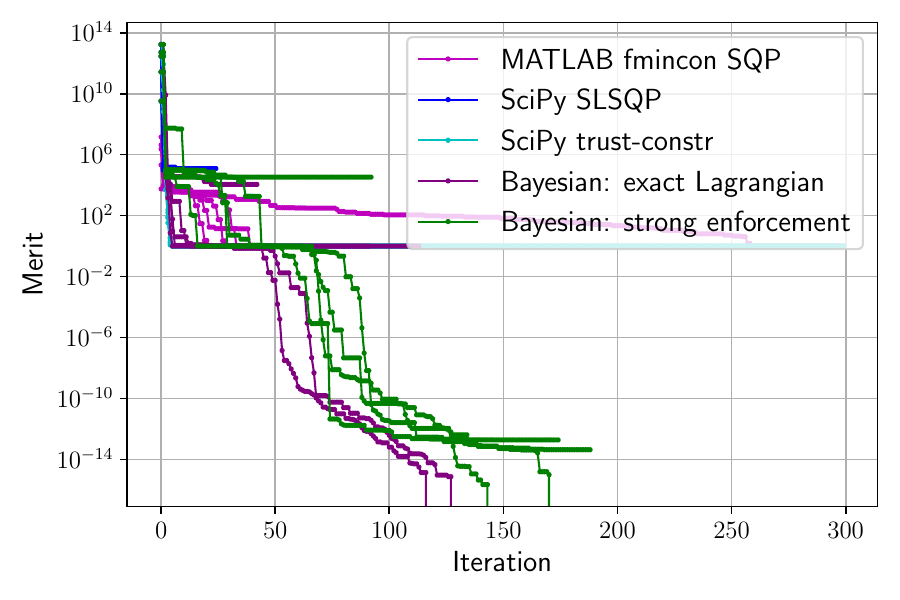}
		\captionsetup{justification=raggedright,singlelinecheck=false}
		\caption{\Eq{Eq_min_prod_nlc_h} with $n_d = 20$: min $f_{\text{prod}}(\xvec)$ with $\xvec^\top \xvec = 1$.}
		\label{Fig_ConOptz_qN_vs_BO_fprod_nlc_h_d20}
	\end{subfigure}	
	\begin{subfigure}[t]{0.328\textwidth}
		\includegraphics[width=\textwidth]{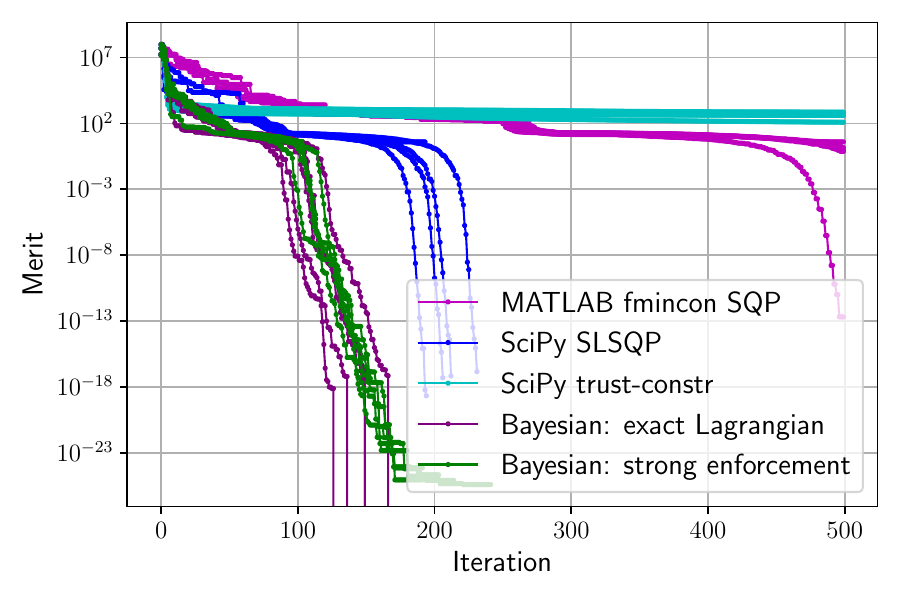}
		\captionsetup{justification=raggedright,singlelinecheck=false}
		\caption{\Eq{Eq_min_nlc_g_Rosen} with $n_d = 20$: min $f_{\text{Rosen}}(\xvec)$ with $\xvec^\top \xvec \leq n_d$.}
		\label{Fig_ConOptz_qN_vs_BO_fRosen_nlc_g_d20}
	\end{subfigure}	
	\caption[Comparison of the Bayesian optimizer with quasi-Newton optimizers for nonlinearly constrained optimization.]{Comparing the Bayesian optimizer to popular quasi-Newton optimizers for nonlinearly-constrained optimization. The merit function is \Eq{Eq_Lag_exa_aug} with $\rho = 100$.}
	\label{Fig_ConOptz_qN_vs_BO}
\end{figure}

Quasi-Newton optimizers are well suited to minimize quadratic functions since they approximate the objective with a quadratic model. Quasi-Newton optimizers generally linearize the nonlinear constraints \cite{nocedal_numerical_2006}, while a separate GP is used to approximate each nonlinear constraint for the Bayesian optimizer. The minimization of the constrained quadratic function $f_{\text{quad}}(\xvec)$ from \Eq{Eq_min_quad_nlc_g} is shown in \Figs{Fig_ConOptz_qN_vs_BO_fquad_nlc_g_d5}{Fig_ConOptz_qN_vs_BO_fquad_nlc_g_d20} for $n_d=5$ and $n_d=20$, respectively. For the $n_d=5$ test case, the SciPy SLSQP optimizer converges the merit function in the fewest number of iterations. The performance of the Bayesian optimizer with both variations is comparable to the SciPy trust-constr and MATLAB fmincon SQP optimizers. However, for $n_d=20$ both variations of the Bayesian optimizer outperform the quasi-Newton optimizers from SciPy and MATLAB. 

The constrained minimization of the product function $f_{\text{prod}}(\xvec)$ from \Eq{Eq_min_prod_nlc_h} is shown for $n_d=5$ and $n_d=20$ in \Figs{Fig_ConOptz_qN_vs_BO_fprod_nlc_h_d5}{Fig_ConOptz_qN_vs_BO_fprod_nlc_h_d20}, respectively. For the $n_d=5$ case, there is a greater variation in the number of iterations required for the SciPy and MATLAB optimizers to reduce the merit function below $10^{-5}$ than the Bayesian optimizer. Even so, these quasi-Newton optimizers take on average about as many iterations to achieve this desired tolerance as the two variations of the Bayesian optimizer. However, for the $n_d=20$ case the quasi-Newton optimizers do not reduce the merit function below $10^{-5}$ for a single optimization run, while the Bayesian optimizer does so for the majority of the optimization runs.

\Figs{Fig_ConOptz_qN_vs_BO_fRosen_nlc_g_d5}{Fig_ConOptz_qN_vs_BO_fRosen_nlc_g_d20} show the constrained minimization of the constrained Rosenbrock function from \Eq{Eq_min_prod_nlc_h} with $n_d=5$ and $n_d=20$, respectively. For both test cases, the two variations of the Bayesian optimizer are able to reduce the merit function below $10^{-5}$ in significantly fewer iterations than the SciPy and MATLAB optimizers. The SciPy trust-constr and MATLAB fmincon SQP optimizers both struggle for this case, especially for $n_d=20$.

%--------------------------------------------------------------------------
%-------------------------------------------------------------------------- 

\begin{figure}[t!]
	\centering
	\begin{subfigure}[t]{0.328\textwidth}
		\includegraphics[width=\textwidth]{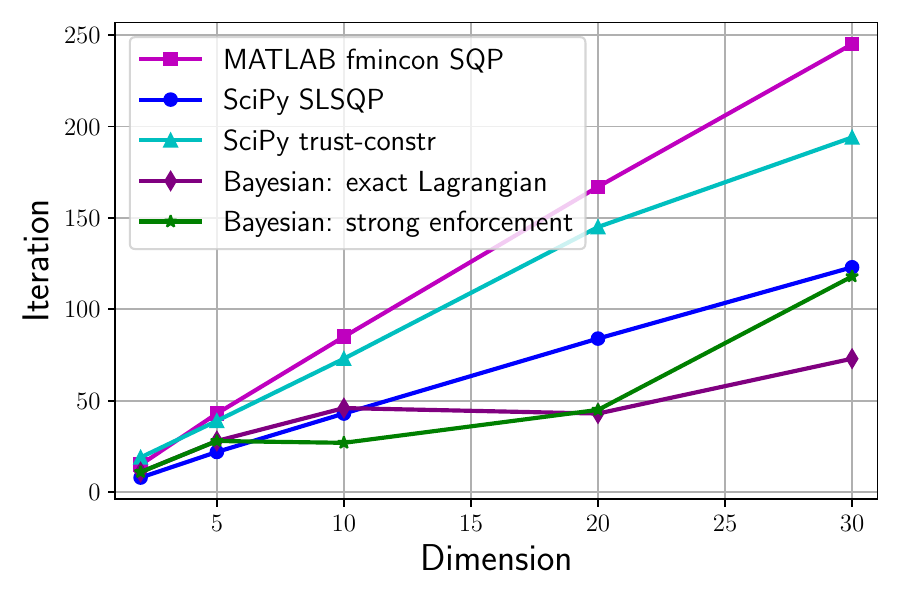}
		\captionsetup{justification=raggedright,singlelinecheck=false}
		\caption{min $f_{\text{quad}}(\xvec)$ with $\xvec^\top \xvec \geq 4$ from \Eq{Eq_min_quad_nlc_g}.}
		\label{Fig_ConOptz_qN_vs_BO_tol_fquad}
	\end{subfigure}	
	\begin{subfigure}[t]{0.328\textwidth}
		\includegraphics[width=\textwidth]{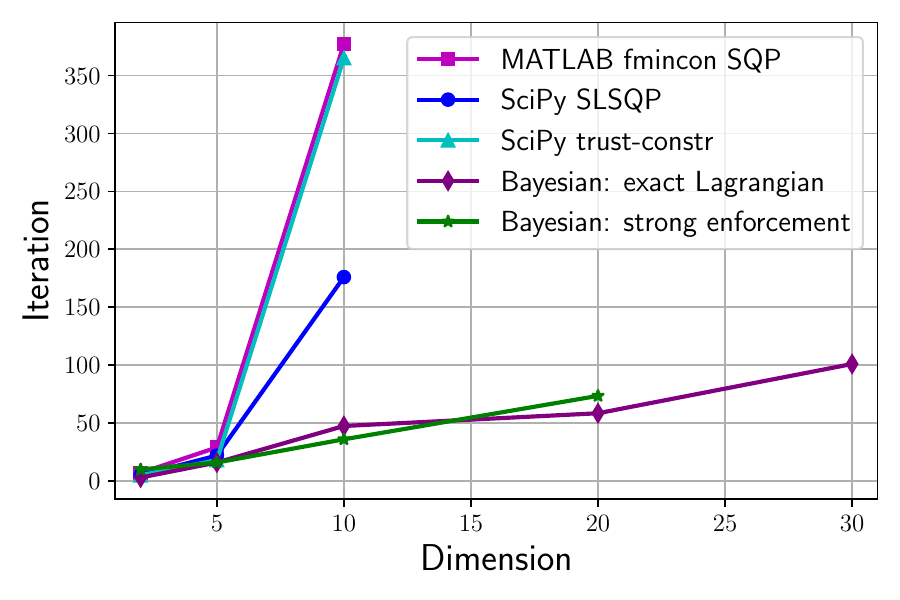}
		\captionsetup{justification=raggedright,singlelinecheck=false}
		\caption{min $f_{\text{prod}}(\xvec)$ with $\xvec^\top \xvec = 1$ from \Eq{Eq_min_prod_nlc_h}.}
		\label{Fig_ConOptz_qN_vs_BO_tol_fprod}
	\end{subfigure}	
	\begin{subfigure}[t]{0.328\textwidth}
		\includegraphics[width=\textwidth]{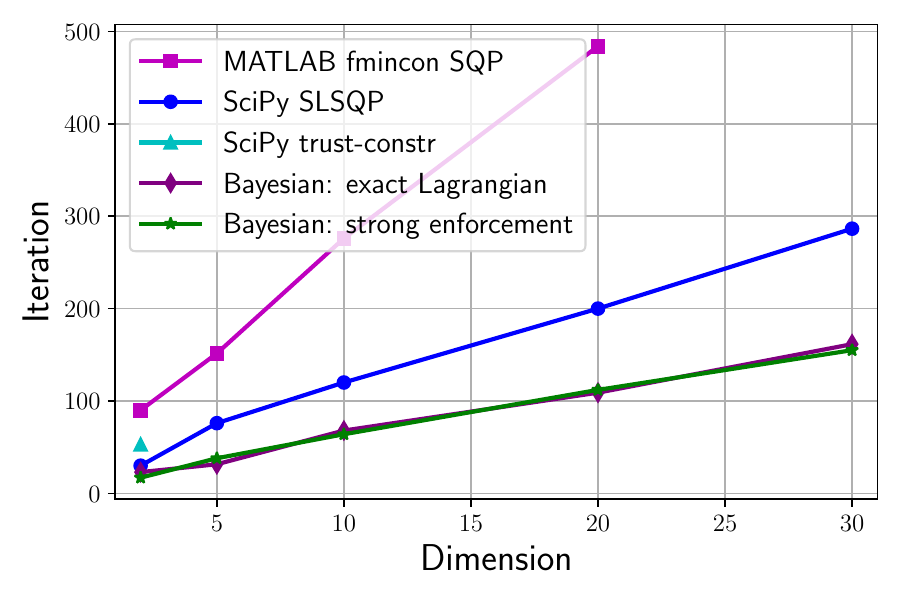}
		\captionsetup{justification=raggedright,singlelinecheck=false}
		\caption{min $f_{\text{Rosen}}(\xvec)$ with $\xvec^\top \xvec \leq n_d$ from \Eq{Eq_min_nlc_g_Rosen}.}
		\label{Fig_ConOptz_qN_vs_BO_tol_fRosen}
	\end{subfigure}	
	\caption[Median number of iterations for the Bayesian optimizer using a strong enforcement of the nonlinear constraints to reduce the merit function below $10^{-5}$.]{The median number of iterations for the Bayesian optimizer using Algorithm~\ref{Alg_StrongEnf} to strongly enforcement the nonlinear constraints to reduce the merit function below $10^{-5}$. The merit function comes from \Eq{Eq_Lag_exa_aug} with $\rho = 100$.}
	\label{Fig_ConOptz_qN_vs_BO_tol}
\end{figure}

The median number of iterations required for the optimizers to reduce the merit function below $10^{-5}$ is plotted in \Fig{Fig_ConOptz_qN_vs_BO_tol}. The Bayesian optimizer with the exact augmented Lagrangian and the strong enforcement of the nonlinear constraints outperforms the SciPy and MATLAB optimizers for all three test cases. This even applies in \Fig{Fig_ConOptz_qN_vs_BO_tol_fquad} to the minimization of the constrained quadratic function $f_{\text{quad}}(\xvec)$ from \Eq{Eq_min_quad_nlc_g}, which the quasi-Newton methods are particularly well suited to minimize. For the constrained minimization of the product function $f_{\text{quad}}(\xvec)$ from \Eq{Eq_min_prod_nlc_h}, which is shown in \Fig{Fig_ConOptz_qN_vs_BO_tol_fprod}, the SciPy and MATLAB optimizers are only competitive with the Bayesian optimizer when $n_d=2$ and $n_d=5$. For the higher-dimensional cases, the quasi-Newton optimizers either take significantly more iterations to achieve the desired tolerance, or are not able to do so for any of the optimization runs. 

%--------------------------------------------------------------------------
%-------------------------------------------------------------------------- 

All of the optimizers performed five independent optimization runs for each of the test cases. \Table{Table_optz_tol_cstr_OptzAll} tabulates how many of these runs achieved the desired tolerance of a merit function smaller than $10^{-5}$. In each cell there are three numbers, which correspond to the constrained test cases from \Eqss{Eq_min_quad_nlc_g}{Eq_min_prod_nlc_h}{Eq_min_nlc_g_Rosen}, respectively. From the first number in each cell it is clear that all of the optimizers were generally able to achieve the desired tolerance for the constrained minimization of the quadratic function $f_{\text{quad}}(\xvec)$. However, the low second number in each cell indicates that all of the optimizers struggled with the constrained minimization of the product function $f_{\text{prod}}(\xvec)$ for $n_d \geq 10$. The SciPy and MATLAB optimizer were not able to achieve the desired tolerance for a single run when $n_d=20$ or $n_d=30$, and only once when $n_d=10$. Finally, the last number in each cell shows that the Bayesian optimizer achieved the desired tolerance more often than the SciPy and MATLAB optimizer for the constrained minimization of the Rosenbrock function.

%--------------------------------------------------------------------------
\begin{table}[t]
	\captionsetup{font=scriptsize}
	\centering
	\begin{tabular}{l ccccc | c} 
		% \hline
		Optimizers & $n_d = 2$ & $n_d = 5$ & $n_d = 10$ & $n_d = 20$ & $n_d = 30$ & Total (out of 25) \\ % \textcolor{Green}{5}
		\hline 
		MATLAB fmincon SQP 			& $\cFive \cdot \cFive \cdot \cFive$ & $\cFive \cdot \cFive \cdot 4$ 		& $\cFive \cdot 1 \cdot \cFive$ & $\cFive \cdot \cZero \cdot 1$ & $\cFive \cdot \cZero \cdot \cZero$ & $ \cTentyFive \cdot 11 \cdot 15$ \\
		SciPy SLSQP 				& $\cFive \cdot \cFive \cdot \cFive$ & $\cFive \cdot \cFive \cdot \cFive$ 	& $\cFive \cdot 1 \cdot \cFive$ & $\cFive \cdot \cZero \cdot 4$ & $\cFive \cdot \cZero \cdot 2$ & $ \cTentyFive \cdot 11 \cdot 21$ \\
		SciPy trust-constr 			& $\cFive \cdot \cFive \cdot 3$ 	 & $\cFive \cdot \cFive \cdot \cZero$ 	& $\cFive \cdot 1 \cdot \cZero$ & $\cFive \cdot \cZero \cdot 0$ & $\cFive \cdot \cZero \cdot \cZero$ & $ \cTentyFive \cdot 11 \cdot 3 \phantom{0}$ \\
		Bayesian: exact Lagrangian 	& $\cFive \cdot \cFive \cdot \cFive$ & $\cFive \cdot \cFive \cdot 4$ 		& $\cFive \cdot 2 \cdot 3 $ 	& $\cFive \cdot 2 \cdot \cFive$ & $1 \cdot 2 \cdot 4$ & $ 21 \cdot 16 \cdot 21$ \\
		Bayesian: strong 			& $\cFive \cdot \cFive \cdot \cFive$ & $\cFive \cdot \cFive \cdot \cFive$ 	& $\cFive \cdot 2 \cdot \cFive$ & $\cFive \cdot 4 \cdot \cFive$ & $\cFive \cdot \cZero \cdot \cFive$ & $ \cTentyFive \cdot 16 \cdot \cTentyFive$ \\
		\hline
	\end{tabular}
	\caption[Number of optimization runs for the quasi-Newton optimizers from SciPy and MATLAB, and the Bayesian optimizer to achive the desired tolerance for the nonlinearly-constrained optimization test cases.]{The number of optimization runs out of five that the SciPy, MATLAB, and Bayesian optimizer reduce the merit function below $10^{-5}$ for the nonlinearly constrained test cases. The set of three numbers in each cell is for the constrained test cases from \Eqss{Eq_min_quad_nlc_g}{Eq_min_prod_nlc_h}{Eq_min_nlc_g_Rosen}, respectively.}
	\label{Table_optz_tol_cstr_OptzAll}
\end{table}
%--------------------------------------------------------------------------

The results in this section demonstrate that the Bayesian optimizer using either the exact augmented Lagrangian or the strong enforcement of the nonlinear constraints is able to reduce the merit function below a desired tolerance in significantly fewer iterations than popular quasi-Newton optimizers from SciPy and MATLAB. Furthermore, the Bayesian optimizer has proven itself to be at least as robust as the SciPy and MATLAB optimizers for the constrained minimization of the quadratic and Rosenbrock functions from \Eqs{Eq_min_quad_nlc_g}{Eq_min_nlc_g_Rosen}, respectively. For the constrained product function from \Eq{Eq_min_prod_nlc_h}, the Bayesian optimizer was significantly more robust than the SciPy and MATLAB optimizers, which did not achieve the desired tolerance once for $n_d \geq 20$.

A simplified version of the strong enforcement method from this paper was used to perform aerodynamic shape optimization of an airfoil at transonic speed. The objective was drag minimization and there were two nonlinear constraints: a target lift and a required minimum area \cite{marchildon_gradient-enhanced_2024}. The Reynolds-averaged Navier-Stokes equations were solved, and the Bayesian optimizer was compared to the quasi-Newton optimizer SNOPT. Both optimizers exhibited similar performance, achieving a five order of magnitude reduction in the optimality and a feasibility below $10^{-8}$.

%% file: Sections/Sec_Conclusions.tex
%!TEX root = ../Marchildon_ZinggCstrBayeOptz.tex

% ------------------------------------------
% New section
% ------------------------------------------
\section{Conclusions} \label{Sec_Conclusions}

This paper was focused on developing a framework to enable efficient local nonlinearly-constrained gradient-enhanced Bayesian optimization. A previously developed unconstrained local Bayesian optimization framework was leveraged and extended to handle nonlinear constraints. This local optimization framework uses a data region to select a subset of evaluation points to train and evaluate the GP and includes a trust region that leverages the probabilistic component of the GP's posterior \cite{marchildon_efficient_2025}. 

Exact augmented Lagrangian and strong enforcement methods were presented that enable effective nonlinearly-constrained Bayesian optimization. The first method uses an acquisition function based on an exact augmented Lagrangian and a probabilistic penalty that promotes exploration when the constraints are not satisfied. The second method modifies the minimization of the acquisition function to include additional constraints. These additional constraints are formed from the posterior of GPs approximating the nonlinear constraints and their bounds get smaller as the feasibility is reduced. Both of these new methods can handle nonlinear inequality and equality constraints.

The new methods were compared to the simple $\ell_2$ penalty method and the popular probability of feasibility methods with both the upper confidence and expected improvement acquisition functions. For the three test cases that were considered with 2 to 30 variables, the exact augmented Lagrangian and the strong enforcement methods were both able to converge the merit function more quickly and to achieve a final tolerance that was several orders of magnitude lower than the $\ell_2$ penalty and probability of feasibility methods. Furthermore, the exact augmented Lagrangian and the strong enforcement methods can be applied to problems containing nonlinear equality constraints, unlike the probability of feasibility method and most other methods developed for nonlinearly-constrained Bayesian optimizers.

The Bayesian optimizer with either the exact augmented Lagrangian method or the strong enforcement method was compared to quasi-Newton optimizers from SciPy and MATLAB for the same set of three nonlinearly-constrained test cases. With either method, the Bayesian optimizer required fewer fewer function evaluations to reach the optimization tolerance of a merit function smaller than $10^{-5}$ and did so more consistently than the quasi-Newton optimizers. These results demonstrate that the local optimization framework that was leveraged along with either of the two new constrained methods enable the Bayesian optimizer to be an effective and robust local minimizer for nonlinearly-constrained problems. 

The characteristic of the optimization problem should be considered carefully when deciding which of the two methods to try first. For a problem where optimizers struggle to converge, the strong enforcement method is a good candidate since it has several parameters that can be tuned intuitively. Alternatively, if a user is interested in solving a multimodal nonlinearly-constrained optimization problem, then the augmented Lagrangian method could be considered since it would be easier to extend to this class of problem. This would require changes to the local optimization framework, such as removing or loosening the trust regions.

%% file: Sections/Appendix_ExaLagMul.tex
%!TEX root = ../Marchildon_ZinggCstrBayeOptz.tex

% ------------------------------------------
% New section
% ------------------------------------------
\section{Closed form solution for Lagrange multipliers} \label{Sec_Appendix_ExaLagMul}

In this appendix we present the closed form solution for the Lagrange multipliers for the exact Augmented Lagrangian from Di Pillo \cite{spedicato_exact_1994}. The Lagrange multipliers are selected by minimizing the function $\Psi(\psivec_h, \psivec_g; \xvec)$ from \Eq{Eq_exaLagMul}, which is quadratic and it thus has a closed-form solution. The solution uses the inverse of the following symmetric positive definite matrix:
\begin{equation} \label{Eq_exaLag_matrixM}
	\M(\xvec) = 
	\begin{bmatrix}
		\M_{11}(\xvec) & \M_{12}(\xvec) \\
		\M_{21}(\xvec) & \M_{22}(\xvec)
	\end{bmatrix},
\end{equation}
which has $\nnlcg + \nnlch$ rows and columns and
\begin{align}
	\M_{11}(\xvec) 
	&= \nabla_x \gvec(\xvec) \left( \nabla_x \gvec(\xvec) \right)^\top + \alpha_1 \G^2(\xvec) + \alpha_2 w(\xvec) \Idty_{\nnlcg} \\
	\M_{12}^\top (\xvec) = \M_{21}(\xvec) 
	&= \nabla_x \gvec(\xvec) \left( \nabla_x \hvec(\xvec) \right)^\top \\
	\M_{22}(\xvec) 
	&= \nabla_x \hvec(\xvec) \left( \nabla_x \hvec(\xvec) \right)^\top + \alpha_2 w(\xvec) \Idty_{\nnlch}.
\end{align}
The Lagrange multipliers are thus calculated at each point $\xvec$ in the parameter space with
\begin{align}
	\begin{bmatrix}
		\psivec_g (\xvec) \\
		\psivec_h (\xvec)
	\end{bmatrix}
	&= \argmin_{\hat{\psivec}_g, \hat{\psivec}_h} \Psi(\hat{\psivec}_g, \hat{\psivec}_h; \xvec) \nonumber \\
	&= - \M^{-1}(\xvec) 
	\begin{bmatrix}
		\nabla_x \gvec(\xvec) \\
		\nabla_x \hvec(\xvec)
	\end{bmatrix}
	\left(\nabla_x f(\xvec) \right)^\top, \label{Eq_ExaLagMul_sol}
\end{align}
where this can be solved efficiently by calculating the Cholesky decomposition of $\M$. If the objective $f(\xvec)$, and nonlinear constraints $\gvec(\xvec)$ and $\hvec(\xvec)$ are all at least twice continuously differentiable, then the derivative of the exact augmented Lagrangian from \Eq{Eq_Lag_exa_aug} is also continuous.

%% file: Sections/Appendix_TestCaseSol.tex
%!TEX root = ../Marchildon_ZinggCstrBayeOptz.tex

% ------------------------------------------
% New section
% ------------------------------------------
\section{Solutions to the constrained minimization problems} \label{Sec_Appendix_TestCase}

% ------------------------------------------
% New subsection
% ------------------------------------------
\subsection{Solution to the constrained minimization of the quadratic function} \label{Sec_Appendix_TestCase_SolNlcQuad}

In this section, the solution to the constrained optimization test case from \Eq{Eq_min_quad_nlc_g} is derived. In \Sec{Sec_ConOptz_CstrTestCases} the test case has a constraint with a radius of 2 while it is generalized in this appendix to be simply $r > 0$:
\begin{equation} \label{Eq_min_quad_A_with_h}
	\min_{\xvec} \xvec^\top \A \xvec - r^2 \eigmin \quad \text{subject to} \quad \| \xvec \|_2^2 = r^2,
\end{equation}
where $\A$ is a symmetric positive definite matrix and $\eigmin$ its smallest eigenvalue. The Lagrangian of this constrained problem is
\begin{equation}
	J(\xvec) = \xvec^\top \A \xvec - r^2 \eigmin + \psi \left( r^2 - \| \xvec \|_2^2 \right).
\end{equation}
The derivative of the Lagrangian with respect to $\xvec$ is 
\begin{align*}
	\p{J}{\xvec}
	&= 2 \A \xvec - 2 \psi \xvec \\
	&= 2 \left( \A - \Idty \psi \right) \xvec, \yesnumber \label{Eq_KKT_convexA}
\end{align*}
which requires $\psi$ and $\xvec$ to be an eigenvalue and an eigenvector of $\A$, respectively. The trivial solution $\xvec = \zero$ cannot be used since it does not satisfy the nonlinear constraint. The objective function can now be evaluated
\begin{align*}
	f(\xvec) 
	&= \xvec^\top \A \xvec - r^2 \eigmin \\
	&= \xvec^\top \left( \psi \xvec \right) - r^2 \eigmin \\
	&= r^2 (\psi - \eigmin) \\
	&=0, %\yesnumber \label{Eq_eval_f_with_psi}
\end{align*}
where $\xvec^\top \xvec = r^2$ was used since it is required for the nonlinear constraint to be satisfied and the objective function is minimized with $\psi = \eigmin$, since $\psi$ must be an eigenvalue of $\A$. The solution to \Eq{Eq_min_quad_A_with_h} thus evaluates to zero at $\xvec = \pm r \uvec_{\min}$, where $\uvec_{\min}$ is the unit eigenvector of $\A$ associated with the eigenvalue $\eigmin$. There are two equivalent solutions since the objective function and nonlinear constraint are unchanged from a reflection through the origin, \ie $f(\xvec) = f(-\xvec)$ and $h(\xvec) = h(-\xvec)$. The solution to \Eq{Eq_min_quad_A_with_h} is the same if the constraint is replaced with $\| \xvec \|_2^2 \geq r^2$ since this inequality constraint would be active at the solution.

% ------------------------------------------
% New subsection
% ------------------------------------------
\subsection{Solution to the constrained minimization of the product function}
\label{Sec_Appendix_TestCase_SolNlcProd}
%\label{Sec_AppendixUnconStudy_SolNlcProd}

We consider the solution to the minimization or the product function from \Eq{Eq_prod_fun} with the solution constrained to be on the circle of radius $r > 0$:
\begin{equation} \label{Eq_min_prod_w_h}
	\min_{\xvec} 1 - n_d^{n_d / 2} \prod_{i=1}^{n_d} x_i \quad \text{subject to} \quad \| \xvec \|_2^2 - r^2 = 0.
\end{equation}
The Lagrangian for this function is 
\begin{equation}
	J(\xvec) = 1 - n_d^{n_d / 2} \prod_{i=1}^{n_d} x_i + \psi \left( \| \xvec \|_2^2 - r^2 \right),
\end{equation} 
which is minimized when its derivatives with respect to the entries of $\xvec$ are zero:
\begin{align}
	\p{J}{x_i} 
	&= -n_d^{n_d / 2} + 2 \psi x_i = 0 \\
	x_i	&= \frac{n_d^{n_d / 2}}{2 \psi} \, \forall \, i \in \{1, \ldots, n_d\},
\end{align}
where the Lagrangian is clearly minimized when all of the entries in $\xvec$ are equal, \ie $\xvec = \alpha \one$. To satisfy the nonlinear constraint we must have
\begin{align*}
	\xvec^\top \xvec 
	&= \alpha^2 \one^\top \one \\
	&= \alpha^2 n_d \\
	&= r^2,
\end{align*}
which is satisfied with $\alpha = \frac{r}{\sqrt{n_d}}$. The solution to \Eq{Eq_min_prod_w_h} is thus $\xvec = \frac{r}{\sqrt{n_d}} \one$, where the objective evaluates to $1 - r^{n_d}$.